\mag=\magstep1
\documentclass[10pt,a4paper]{amsart}
\usepackage{latexsym,amssymb,amscd}

\def\over{/}
\def\Cal{\mathcal}
\def\frak{\mathfrak}

\def\<<{\langle}
\def\>>{\rangle}
\def\Diag {\operatorname{Diag}}

\numberwithin{equation}{section}
\newtheorem{theorem}{Theorem}[section]
\newtheorem{proposition}[theorem]{Proposition}
\newtheorem{corollary}[theorem]{Corollary}
\newtheorem{definition}[theorem]{Definition}

\newtheorem{remark}[theorem]{Remark}
\newtheorem{lemma}[theorem]{Lemma}

\newtheorem{notation}[theorem]{Notation}

\begin{document}

\title{
Theta constants associated to coverings of $\mathbf P^1$
branching at $8$ points
}

\author{Keiji Matsumoto}
\address{Division of Mathematics, Graduate School of Science, 
Hokkaido University, Japan}
\email{matsu@math.sci.hokudai.ac.jp}

\author{Tomohide Terasoma}
\address{Department of Mathematical Science, University of Tokyo,
Komaba, Meguro, Japan}
\email{terasoma@ms.u-tokyo.ac.jp}

\subjclass{
Primary  14J15; 
Secondary 32N15,11F55
}
\date{\today }
\keywords{period map, theta constant, Prym variety}

\maketitle

\makeatletter
\renewcommand{\@evenhead}{\tiny \thepage \hfill  K.MATSUMOTO and T.TERASOMA 
\hfill}
\renewcommand{\@oddhead}{\tiny \hfill  THETA FUNCTION OF BRANCHED COVER
 \hfill \thepage}
\makeatother

\begin{abstract}
In this paper, we construct automorphic forms on the five dimensional
complex ball which give the inverse of the period map
for cyclic 4-ple coverings of the complex projective line
branching at eight points. 
We use theta constants associated to the Prym varieties 
of these coverings.
\end{abstract}

\section{Introduction}
Let $\mu_1,\dots,\mu_n$ be rational numbers such that $0<\mu_j<1$, 
$\sum_{j=1}^n\mu_j=2$, and 
let $d$ be the common denominator of $\mu_1,\dots,\mu_n$. 
For the cyclic $d$-ple covering $C$ of $\mathbf P^1$ 
branching at $n$ points with index $\mu=(\mu_1,\dots,\mu_n)$ and 
a homology marking $\phi$ of  $C$, 
the period $p(C,\phi)$ of a marked curve $(C,\phi)$ 
can be regarded as an element of the $(n-3)$-dimensional complex ball $B_\mu$.
The morphism $p$ from the moduli space $M_{marked}$ of marked curves 
$(C,\phi)$ to $B_\mu$ is called a period map. 
Since the period map $p$ is equivariant under the monodromy group 
$\Gamma_\mu$, it induces the morphism 
$M_{marked}/\Gamma_\mu\to B_\mu/\Gamma_\mu$. 
According to results of \cite{DM} and \cite{T}, the period map $p$ is 
an isomorphism onto a Zariski open set of $B_\mu/\Gamma_\mu$ 
if $\mu$ satisfies the condition 
\begin{equation}
\label{DMT}
(1-\mu_j -\mu_k)^{-1} \in \mathbf Z\cup {\infty}\text{ for $j\neq k$ }. 
\end{equation}
There are finitely many such $\mu$'s for $n\ge 5$; 
for $n\ge 9$ there is no $\mu$, 
for $n=8$ there is one,
for $n=7$ there is one,
for $n=6$ there are seven and 
for $n=5$ there are 27.

It is a natural demand to describe the inverse of the period map in terms of 
explicit automorphic forms with respect to $\Gamma_\mu$ for such $\mu$. 
In fact, for several $\mu$ for $n=5,6,$ the inverse of the period map 
is expressed in terms of theta constants; refer to  \cite{K1},  
\cite{Ma1}, \cite{Ma2} and \cite{S}.  
In this paper, we construct automorphic forms on $B_{\mu}$
which give the inverse of the period map for 
the case $n=8,\ \mu=(1/4)^8=({1\over 4},\dots,{1\over 4})$.  
Note that 
all $\mu$'s with $d=4$ satisfying the condition (\ref{DMT}) 
can be obtained by confluences of the branching index $(1/4)^8$.

Before stating our main theorem,
we define a projective embedding of the moduli space $M_{8pts}$
of $8$ points on $\mathbf P^1$, which is isomorphic to 
$M_{marked}/\Gamma_\mu$. 
We divide the set $\{1,\dots,8\}$ into a set of four pairs 
$\{\{j_1, j_2\},\dots \{j_7, j_8\}\}$,  
which is called a $(2,2,2,2)$-partition of $\{1,\dots,8\}$.   
We denote the set of $(2,2,2,2)$-partitions of $\{1,\dots,8\}$
by $P(2^4)$ which has cardinality $105$.
We associate a polynomial 
$P_{r}=\prod_{k=1}^4(x_{j_{2k-1}}-x_{j_{2k}})$ 
for each $(2,2,2,2)$-partition $r=\{\{j_1, j_2\},\dots \{j_7, j_8\}\}$.  
If we regard $x_j$'s as affine coordinates of eight points on $\mathbf P^1$,  
then $P_{r}$ are relative invariants under projective transformations of 
$\mathbf P^1$. 
Thus the set of polynomials
$\{P_r\}_{r \in P(2^4)}$
induces a map $P:M_{8pts} \to \mathbf P^{104}$. 
It is shown in \cite{K2} that  $P$ is an embedding.  

\begin{theorem}[Main Theorem]
There exist $105$ automorphic forms $\Cal T_r^{(2)}$ 
with respect to $\Gamma_\mu$ such that 
the following diagram commutes 

\setlength{\unitlength}{0.75mm}
\begin{picture}(200,45)(0,-5)
\put(34,0){$B_{\mu}/\Gamma_\mu$}
\put(38,30){$M_{8pts}$}
\put(38,14.5){$p$}
\put(71,14.5){$\mathbf  P^{104}$,}
\put(55,27){$P$}
\put(58,4){$\Cal T^{(2)}$}
\put(42,25){\vector(0,-1){20}}
\put(48,4){\vector(2,1){20}}
\put(48,28){\vector(2,-1){20}}
\end{picture}
where the map $\Cal T^{(2)}$ 
is given by the $105$ automorphic forms $\Cal T_r^{(2)}$.    
Especially, the image of $\Cal T^{(2)}$ coincides with that of $P$ and 
$\Cal T^{(2)}\circ P^{-1}$ gives the inverse of $p$.  
\end{theorem}

Let us explain how to construct $105$ automorphic forms 
in terms of theta constants. 
The curve  $C$ is of genus $9$ and it can be regarded as 
a double cover of a hyperelliptic curve of genus $3$.  
The period $p(C,\phi)$ is an element of the $5$-dimensional complex ball.  
We consider the Prym variety $Prym(C)$ of $C$ with respect to $\rho^2$, 
which is a $6$-dimensional sub-Abelian variety of the Jacobian $J(C)$ 
of $C$ obtained by 
the $(-1)$-eigenspaces of $H^0(C, \Omega^1)$ and $H_1(C, \mathbf  Z)$ 
for the action of $\rho^2$, where $\rho$ is a generator of the group of 
covering transformations of $C \to \mathbf P^1$. 

Since the polarization of $Prym(C)$ is not principal,   
we construct $105$ principally polarized abelian varieties isogenous to 
$Prym(C)$ as follows. 
Let  $Prym(C)_{1-\rho}$ be the group of $(1-\rho)$-torsion points 
of $Prym(C)$, which is isomorphic to $\mathbf F_2^6$ with a quadratic form. 
There are $105$ three dimensional totally singular subspaces 
$\Lambda_r$ of $Prym(C)_{1-\rho}$. 
For each $\Lambda_r$, $A_{L_r}=Prym(C)/\Lambda_r$ is principally polarized.
We show that there is an isomorphism between the automorphism groups 
$\frak S_8$ of the marking of $8$ points and the orthogonal group 
$O_6^+(\mathbf F_2)$ of
$Prym(C)_{1-\rho}$ yields a one to one correspondence 
between $P(2^4)$ and the set $\{\Lambda_r\}$.

For $\Lambda_1$ corresponding to $r_1=\{\{12\},\{34\},\{56\},\{78\}\}
\in P(2^4),$ we study the behavior of 
the pull back $F_{m_1}=\iota^*(\vartheta_{m_1})$ of a theta function 
$\vartheta_{m_1}(z)$ associated to $A_{L_1}$ with a characteristic $m_1$ 
under the composite map 
$\iota:C\to A_{L_1}$ of the canonical map $jac^-:C\to Prym(C)$ and 
the natural projection $Prym(C)\to A_{L_1}$ 
(see Section \ref{Theta principal}).
There are twelve zeros of $F_{m_1}$ in $C$, 
eight of them are in the set of fixed points of $\rho$.
There are three theta functions $\vartheta_{m_j}(z)$ $(j=2,3,4)$ on $A_{L_1}$ 
such that the order of zero of $F_{m_j}=\iota^*(\vartheta_{m_j})$ 
at every fixed point of $\rho$ is same as that of $F_{m_1}$. 
By considering the four zeros of $F_{m_j}$ not fixed by $\rho$, 
we can express the cross ratio of $x_1,x_2,x_5,x_6$ in terms of 
the theta constants $\vartheta_{m_1},\dots,\vartheta_{m_4}$.

The product $\Cal T^{(2)}_1=\prod_{j=1}^4\vartheta_{m_j}$ 
is an automorphic form with respect to the monodromy group  $\Gamma_\mu$.
The period map $p$ induces an isomorphism of groups $Aut(M_{8pts})\simeq
Aut(B_{\mu}/\Gamma_\mu)$. For each $(2,2,2,2)$-partition $r\in P(2^4)=
Stab(r_1)\backslash \frak S_8$, 
we define $\Cal T_r^{(2)}$ by using the action of  
$\frak S_8\subset Aut(M_{8pts})\simeq Aut(B_{\mu}/\Gamma_\mu)$, 
where $Stab(r_1)$ is the stabilizer of the partition $r_1$. 
In order to prove our main theorem, we investigate 
the action of $Stab(r_1)$ on the space generated by the theta constants
and relations between theta functions for $A_{L_1}$ and those for $A_{L_r}$.

The authors are grateful that Shigeyuki Kondo informed us 
of his research on constructing automorphic forms on the $5$-dimensional 
complex ball in terms of Borcherds products.

\begin{notation}
\label{beginig notation}
In this paper, the imaginary unit is denoted by $i$.
For an element $\alpha \in \bold C$,
$\exp (2\pi i \alpha)$ is denoted by   
$\bold e(\alpha )$.
For a square matrix $A$, the vector consisting of the diagonal elements
of $A$ is denoted by $A_0$. For a vector $v=(v_1, \dots, v_k)$, 
the diagonal matrix 
$
\begin{pmatrix}
v_1 & & \\
&\ddots & \\
& & v_k
\end{pmatrix}
$
is denoted by ${\rm diag}(v)$.
\end{notation}

\section{The Prym variety of $C$}
\subsection{$4$-ple covering of $\mathbf  P^1$ branching at $8$ points}
\label{homology of covering br at 8 pts}
Let $C$ be the projective smooth model of an algebraic curve defined by 
\begin{equation}
\label{curve C}
w^4 = \prod_{j=1}^8 (z-x_j),
\end{equation}
where $x_1, \dots, x_8$ are distinct elements of $\mathbf  C$.
The curve $C$ is of genus $9$  since it can be regarded as a $4$-ple 
covering of $\mathbf P^1$ branching at eight points. The automorphism 
\begin{equation}
\label{equation of branched curve}
\rho:C\ni (z,w) \mapsto (z,iw)\in C
\end{equation}
induces actions on $H^1(C,\mathbf Q)$ and $H_1(C,\mathbf Q).$ 
We denote the $(-1)$-eigenspaces of 
$H^1(C,\mathbf Q)$ and $H_1(C,\mathbf Q)$ of $\rho^2$ by 
$H^1(C,\mathbf Q)^-$ and $H_1(C,\mathbf Q)^-,$ 
respectively. We put 
$$H^1(C, \mathbf  Z)^-=H^1(C, \mathbf  Q)^- \cap H^1(C, \mathbf  Z), \quad 
H_1(C, \mathbf  Z)^-=H_1(C, \mathbf  Q)^- \cap H_1(C, \mathbf  Z).$$
Since the action $\rho$ preserves the
polarized rational Hodge structure of $H^1(C, \mathbf  Q)$, the $(-1)$-eigen
subspace $H^1(C, \mathbf  Q)^-$ of $\rho^2$ is a polarized rational
sub-Hodge structure of $H^1(C, \mathbf  Q)$.
The action of $\rho$ induces an action on each factor of 
the Hodge decomposition
$$
H^1(C, \mathbf  Z)^- \otimes \mathbf C \simeq H^0(C, \Omega^1)^-
\oplus \overline{H^0(C, \Omega^1)^-}.
$$

\begin{proposition}
\label{signature is 1.5}
The multiplicity of the eigenvalue $i$ (resp. $-i$) of 
$\rho$ on $H^0(C, \Omega^1)^-$ is $5$ (resp. $1$). 
\end{proposition}
\begin{proof}
Differential $1$-forms
$$
\varphi_j=\frac{z^j dz}{w^3} \ (j=0, \dots ,4), \quad
\varphi_5=\frac{dz}{w},\quad \varphi_6=\frac{dz}{w^2},\ 
\ \varphi_7=\frac{zdz}{w^2},\ \varphi_8=\frac{z^2dz}{w^2}
$$
span the space $H^0(C, \Omega^1)$. We have
$\rho^*(\varphi_j)=i\varphi_j$ for $j=0, \dots ,4$,
$\rho^*(\varphi_5)=-i\varphi_5$ and 
$\rho^*(\varphi_k)=-\varphi_k$ for $k=6,7,8$. 
\end{proof}

We study the intersection form on 
$H_1(C, \mathbf  Z)^-$. For the moment, we assume
that $x_j\in \mathbf R$ $(j=1,\dots,8)$ and $x_1<x_2<\cdots<x_8.$
Let $\Cal U_0$ be $\mathbf P^1$
cut along the eight semi-lines 
$$
l_j=\{x_j+it \in \mathbf C\subset \mathbf P^1\mid  t \leq 0\}, \ j=1,\dots,8.
$$
The curve $C$ can be regarded as gluing
$\Cal U_0$ and the copies $\Cal U_k=\rho^k(\Cal U_0)$ $(k=1,2,3)$ along 
$l_1, \dots ,l_8$.  
The branch of the function $w$ is assigned 
as its value is in $i^k \mathbf  R_+$ for 
$z \in (-\infty,x_1)\subset \Cal U_k$.
Let $\alpha_j$ ($1 \leq j\leq 7$) be the interval 
$[x_j,x_{j+1}] \subset \mathbf  R$ in the sheet $\Cal U_0$ 
and $\alpha_8$ be 
$[-\infty,x_1]\cup [x_8, \infty]\subset \mathbf  R$ 
in the sheet $\Cal U_0$. 
The orientation of $\alpha_j$ is given in the figure.

\setlength{\unitlength}{0.75mm}
\begin{picture}(200,50)(-20,5)
\put(-3,4){$w \in i^k\cdot \mathbf  R_+$}
\put(5,10){\vector(1,2){8}}
\put(0,26){\dashbox{1}(20,0)}
\put(14,30){$x_1$}
\put(20,6){\dashbox{1}(0,20)}
\put(40,6){\dashbox{1}(0,20)}
\put(60,6){\dashbox{1}(0,20)}
\put(20,26){\vector(1,0){20}}
\put(40,26){\vector(1,0){20}}
\put(57,26){\vector(1,0){3}}
\put(37,26){\vector(1,0){3}}
\put(29,22){$\alpha_1$}
\put(29,10){$\rho\alpha_1$}
\put(39,30){$x_2$}
%
\put(70,46){\vector(-1,0){20}}
\put(50,40){$\alpha_8$}
\put(70,36){\vector(-1,0){20}}

\put(48,30){$\dots$}
\put(59,30){$x_7$}
\put(60,26){\vector(1,0){20}}
\put(77,26){\vector(1,0){3}}
\put(69,22){$\alpha_7$}
\put(67,10){$\rho^7\alpha_7$}
\put(81,30){$x_{8}$}
\put(80,6){\dashbox{1}(0,20)}

\put(50,26){\oval(60,20)[t]}
\put(20,31){\oval(20,30)[l]}
\put(80,31){\oval(20,30)[r]}
\put(20,46){\line(1,0){60}}

\put(20,16){\vector(1,0){20}}
\put(40,16){\vector(1,0){20}}
\put(60,16){\vector(1,0){20}}

\end{picture}
\begin{center}
{\bf Figure 1. }
\end{center}
\vskip 0.2in
Then the $1$-chain $A_j = (1-\rho ^2)\alpha_j$ is a cycle satisfying 
$\rho^2(A_j)=-A_j$. 
We put $B_j = \rho A_j$.
The above figure shows that both of $\sum_{j=1}^8\alpha_j$ 
and $\sum_{j=1}^8\rho^j\alpha_j$ are boundaries of $2$-chains.
Therefore we have
$\sum_{j=1}^8A_j=0$ and $\displaystyle\sum_{j=1}^7\frac{1-\rho^j}{1-\rho}A_j=0$
in $H_1(C, \mathbf  Z)^-$. 
The intersection matrix for $\{ A_j, B_j\}_{j=1, \dots ,6}$ is
given as
\begin{equation}
\label{intersection matrix for A,B}
\left(
\begin{matrix}
P & Q \\
-Q & P
\end{matrix}
\right),
\end{equation}
where
$$
{\small
P= \left(
\begin{matrix}
 0& 1& 0& 0& 0& 0 \\
-1& 0& 1& 0& 0& 0 \\
 0&-1& 0& 1& 0& 0 \\
 0& 0&-1& 0& 1& 0 \\
 0& 0& 0&-1& 0& 1 \\
 0& 0& 0& 0&-1& 0 \\
\end{matrix}
\right),\ 
Q= \left(
\begin{matrix}
 2&-1& 0& 0& 0& 0 \\
-1& 2&-1& 0& 0& 0 \\
 0&-1&2 &-1& 0& 0 \\
 0& 0&-1& 2&-1& 0 \\
 0& 0& 0&-1& 2&-1 \\
 0& 0& 0& 0&-1& 2 \\
\end{matrix}
\right).}
$$

\begin{proposition}
\label{by cutting get base}
The set $\{ A_j, B_j \}_{j=1, \dots ,6}$ is a basis of $H_1(C, \mathbf  Z)^-$.
\end{proposition}
\begin{proof}
The determinant of the matrix (\ref{intersection matrix for A,B}) 
is $2^6$. By the Fey's
result \cite{F}, we have this proposition.
\end{proof}

\subsection{Polarization of the Prym variety}
\label{polarization of prym var}
The polarized Hodge structure of $H^1(C, \mathbf  Z)$ defines the
abelian variety
$$
Prym(C)=Prym(C,\rho^2)=(H^0(C, \Omega^1)^-)^*/H_1(C, \mathbf  Z)^-,
$$
which is called the Prym variety of $C$.
Since the first homology group $H_1(Prym(C), \mathbf  Z)$
of the Prym variety $Prym(C)$ is isomorphic to $H_1(C, \mathbf  Z)^-$,
the restriction 
$$
(,): H_1(C, \mathbf  Z)^- \times H_1(C, \mathbf  Z)^- \to \mathbf  Z
$$
of the intersection form on $H_1(C, \mathbf  Z)$ gives a polarization
of $Prym(C)$. Thus $H_1(C, \mathbf  Z)^- \simeq H_1(Prym(C), \mathbf  Z)$
naturally has a polarized Hodge structure of weight $(-1)$.

\begin{definition}
Let $H$ be a polarized $\mathbf  Z$-Hodge structure of weight $(-1)$.
If $H$ has a basis whose intersection matrix is 
\begin{equation}
\label{normal form of sympl. form}
\left(\begin{matrix} 0& e \\ -e & 0 \end{matrix}\right) ,\quad 
e = {\rm diag}(\varepsilon_1, \dots ,\varepsilon_d),
\end{equation}
then the type of the polarization of $H$ is said to be 
$(\varepsilon_1, \dots ,\varepsilon_d)$.
For a polarized abelian variety $(A, (,))$, the type of
$(H_1(A, \mathbf  Z),(,))$ is called the type of the polarized 
abelian variety.
The polarization of
type $(1, \dots ,1)$ is called principal.
\end{definition}
The curve $C$ can be regarded as a double covering of 
a hyperelliptic curve of genus $3$ branching at eight points. 
The results in \cite{F} implies that 
the type of the polarized abelian variety $Prym(C)$
is  $(2,2,2,1,1,1)$ (see also \cite{Mu}, \S 3 Cor 1). 
We give a symplectic basis $\Sigma$ of 
$H_1(C, \mathbf  Z)^-$.
\begin{proposition}
\label{base for basic principal lattice}
Let $\Sigma = \{ \alpha'_1, \dots ,\alpha'_6, \beta'_1, \dots ,\beta'_6 \}$
be 
\begin{align*}
& \alpha'_1= 2A_1+ 2A_3 + A_4+ B_4, \\
& \alpha'_2= 2A_1+ A_2 + 2A_5 + 2A_6 + B_2 -2B_4 -2B_5, \\
& \alpha'_3=  A_1+ 2A_2+  A_3 + A_5+ 2A_6 -B_1 + B_3 -B_5, \\
& \alpha'_4= A_1, \ 
 \alpha'_5= A_1 + A_3, \ 
 \alpha'_6= -A_1 -A_3 + B_5, \\
& \beta'_1= 2A_1 + 2A_3 + A_5+ 2A_6  -B_5, \\
& \beta'_2= -A_1 -2A_2 -2A_5-2A_6+  B_1 + 2B_4+ 2B_5, \\
& \beta'_3=  A_4 -A_6 + B_4+ 2B_5+ B_6, \\
& \beta'_4=  A_2,  \ 
 \beta'_5=  A_4,   \ 
 \beta'_6=  A_6.    
\end{align*}
Then $\Sigma$ is a basis of $H_1(C, \mathbf  Z)^-$ whose intersection matrix
is given by (\ref{normal form of sympl. form}) for 
$e={\rm diag}(2,2,2,1,1,1)$.
\end{proposition}
To simplify the notation, the half of the intersection form
on $H_1(C, \mathbf  Z)^-$ is denoted by $<,>$.
To classify the sub-principally 
polarized Hodge structures of $(H_1(C, \mathbf  Z),(,))$ of type $(2,2,2,2,2,2)$
is equivalent to classify principally polarized sub-Hodge structures
of $(H_1(C, \mathbf  Z),<,>)$.

For a polarized Hodge structure $H$, the dual Hodge structure
$H^{\perp}$ is defined as
$$
H^{\perp} = \{ v \in H\otimes \mathbf  Q \mid 
<v, w> \in \mathbf  Z \text{ for all } w \in H\}.
$$
It is easy to see that a polarized Hodge structure $H$ is principal
if and only if $H = H^{\perp}$.
From now on, we use the polarization $<,>$ for sub-Hodge structures
of $H_1(C, \mathbf  Z)^-$.
\begin{proposition}

\begin{enumerate}
\item 
$
(H_1(C, \mathbf  Z)^-)^{\perp} = (1-\rho)H_1(C, \mathbf  Z)^-.
$
\item
A principally polarized sub-Hodge structure $L$ of $H_1(C, \mathbf  Z)^-$
contains $(1-\rho)H_1(C, \mathbf  Z)^-$ and is stable under the action
of $\rho$.
\end{enumerate}
\end{proposition}

In the next section, we give a combinatorial description of
the set of principally polarized sub-Hodge structures
of $H_1(C, \mathbf  Z)^-$.

We close this section by giving an example of a principally 
polarized sub-Hodge structure of $H_1(C, \mathbf  Z)^-$
using the basis $A_j, B_j$, ($j=1, \dots , 6$) given in
Proposition \ref{by cutting get base}. 
\begin{proposition}
\label{example good lattice}
The sub-Hodge structure $L_1$ of $H_1(C, \mathbf  Z)^-$ generated by
$(1-\rho)H_1(C, \mathbf  Z)^-$ and $A_1, A_3, A_5$
is principal. Actually the set
$\Sigma_1=\{a_1, \dots ,a_6,b_1, \dots ,b_6\}$, where
\begin{eqnarray*}
a_1&=&A_1,\\
a_2&=&A_1+A_2+B_2,\\
a_3&=&A_1+A_2+B_2+B_3,\\
a_4&=&A_1+A_2-A_4+B_2+B_3+B_4,\\
a_5&=&A_1+A_2+A_5+B_2+B_3,\\
a_6&=&A_1+A_2+A_5+A_6+B_2+B_3+B_6,\\
b_1&=&-B_1,\\
b_2&=&A_2-B_1-B_2,\\
b_3&=&-A_2-A_3-A_4+B_1+B_2-B_4,\\
b_4&=&-A_2-A_3+B_1+B_2,\\
b_5&=&A_2+A_3-B_1-B_2-B_5,\\
b_6&=&A_2+A_3+A_6-B_1-B_2-B_5-B_6,
\end{eqnarray*}
is a symplectic basis of type $(1,1,1,1,1,1)$.
The action $\rho$ on this basis is given by
\begin{equation}
\label{action rho and U}
^t(a_1, \dots,a_6,  b_1, \dots, b_6) 
 \mapsto
\left(\begin{matrix}
0 & -U \\
U & 0 
\end{matrix}\right)
\;^t(a_1, \dots,a_6,  b_1, \dots, b_6),
\end{equation}
where
$$
U = \left(\begin{matrix}
1 & 0 & 0 & 0 & 0 & 0 \\
0 & 1 & 0 & 0 & 0 & 0 \\
0 & 0 & 0 &-1 & 0 & 0 \\
0 & 0 &-1 & 0 & 0 & 0 \\
0 & 0 & 0 & 0 & 1 & 0 \\
0 & 0 & 0 & 0 & 0 & 1 
\end{matrix} \right).
$$
\end{proposition}
\begin{definition}
A pair $(L, \Sigma_L)$ of a principal sub-Hodge structure 
$L$ of $H_1(C, \mathbf  Z)^-$
and a symplectic basis $\Sigma_L$ of $L$ is called a good basis
if the action $\rho$ on $\Sigma_L$ is given by
(\ref{action rho and U}).
\end{definition}
\begin{remark}
\label{torsion}
Using a good base $\Sigma_{L}$, an element of $\frac{1}{1-\rho}L/L$
can be written as $\frac{1}{2}(\mu, \mu U)$ mod $L$ 
($\mu \in \bold Z^6$).
\end{remark}

We will use the lattice $L_1$ and its symplectic basis $\Sigma_1$ in
Proposition \ref{example good lattice} for explicit calculations of 
theta constants.

\section{$O^+_6(2)$-Level structure and $\frak S_8$-marking}
\subsection{Total singular subspaces, length $0$ elements}
\label{total singular}
We define the standard lattice $H_{std}$ as the free
$\mathbf  Z$-module generated by $\bar A_j, \bar B_j$ $(j=1,\dots,6)$.
We introduce an alternating form $<,>$ and an action of $\rho$ 
on $H_{std}$ by the half of the matrix 
(\ref{intersection matrix for A,B})
in \S \ref{homology of covering br at 8 pts}
and $\rho(\bar A_j)= \bar B_j,$ $\rho(\bar B_j)=-\bar A_j$.
Since the action of $\rho^2$ on $H_{std}$
is equal to the multiplication by $(-1)$, $H_{std}$ (resp.
$H_{std, \mathbf  R}=H_{std}\otimes \mathbf  R$) becomes 
a $\mathbf  Z[\rho]/(\rho^2+1)$-module
(resp. a vector space over 
$\mathbf  R(\rho)\simeq \mathbf  C$). 
We define a bilinear form $h(x,y)$ on $H_{std}$ 
by $h(x, y) = <x, \rho y> - <x,y>i$.
Then $h(,)$ becomes a hermitian form
of the signature $(5,1)$ with respect to the complex structure given by
$\mathbf  R(\rho)$.
The group of isomorphisms of $H_{std}$ preserving the alternating form
$<,>$ and the action of $\rho$ is denoted by $U(H_{std})$.
The value of associated Hermitian metric $\tilde q(x)=h(x,x)$ on
$H_1(C, \mathbf  Z)^-$ is integral.
The class $q(x)$ of $\tilde q(x)$ modulo $2$ defines a 
$\mathbf  Z/2\mathbf  Z$-valued quadratic form on 
$H_1(C, \mathbf  Z)^-/(1-\rho)H_1(C, \mathbf  Z)^-$.
This quadratic form $q$ has the following simple form.

Let $(V,q)$ be a $6$-dimensional vector space $V$ 
over $\mathbf  F_2$ with the quadratic form $q$ of Witt defect $0$.
Such $(V,q)$ is constructed as follows. 
The Hamming length of a vector $x \in {\mathbf F}_2^8$ 
is defined by the number of nonzero elements. 
The subset $\tilde V$ consisting of vectors with the even Hamming length 
is a subspace of ${\mathbf F}_2^8$ containing $(1,\dots ,1)$.
Let $V$ be the quotient space of $\tilde V$ by the subspace 
$(1,\dots ,1)\cdot \mathbf  F_2$. 
The half of the Hamming length modulo $2$ becomes 
the quadratic form $q$ on $V$ of Witt defect $0$.

We have an isomorphism 
\begin{equation}
\label{marking mod 2}
(H_1(C, \mathbf  Z)^-/(1-\rho)H_1(C, \mathbf  Z)^-,q) \to (V,q) 
\end{equation}
by mapping the class of $A_j$ to that of $e_j-e_{j+1}$, 
where $e_j$ is the $j$-th unit vector of ${\mathbf F}_2^8$.
The orthogonal group of $(V,q)$ is denoted by $O^+_6(2)$.
The symmetric group $\frak S_8$ of degree $8$ acts on the space $V$ by 
$e_j\cdot \sigma=e_{j\sigma}$ for $\sigma\in \frak S_8$ and 
$j=1, \dots ,8$.
This action defines a group homomorphism $\frak S_8 \to O^+_6(2)$.

\begin{definition}
For  natural numbers $\varepsilon_1, \dots ,\varepsilon_d$ such that 
$\varepsilon_1+\cdots+\varepsilon_d =8$,
a (unordered) partition of $\{1, \dots ,8\}$ 
into sets of cardinality $\varepsilon_1, \dots,\varepsilon_d$
is called an $(\varepsilon_1, \dots ,\varepsilon_d)$-partition. The set of
$(\varepsilon_1, \dots ,\varepsilon_d)$-partitions of $\{1, \dots ,8\}$ is
denoted by $P(\varepsilon_1, \dots ,\varepsilon_d)$. 
The sets $P(2,2,2,2)$ and $P(4,4)$ are
denoted by $P(2^4)$ and $P(4^2)$, respectively.
\end{definition}

Note that $\# P(2^4) =105$ and $\# P(4^2)=35$.
Therefore we have the following propositions.

\begin{proposition}
\begin{enumerate}
\item
The map $\frak S_8 \to O^+_6(2)$ is an isomorphism.
\item
By mapping an element 
$s=\{\{s_1,\dots,s_4\},\{s_5,\dots,s_8\}\} \in P(4^2)$
to a non-zero element $v=e_{s_1}+\cdots +e_{s_4}$ with $q(v)=0$, we have a
one-to-one correspondence
$$
P(4,4) \simeq \{v\in V \mid q(v)=0, v\neq 0\}.
$$
Under this correspondence, the stabilizer of $s\in P(4,4)$ is isomorphic to
the stabilizer of $v$.
\end{enumerate}
\end{proposition}

For an element $I=\{ \{j_1, j_2\}, \dots ,\{j_7,j_8\}\}$ of
$P(2^4)$, we define a subspace 
$V_I=<e_{j_1}-e_{j_2}, \dots ,e_{j_7}-e_{j_8}>$ of $V$.

\begin{proposition}
\label{one to one cor  of index}
Let $\psi : H_1(C, \mathbf  Z)^- \to V$ be the
natural projection. Then $\psi^{-1}(V_I)$ is a principally polarized
sub-Hodge structure in $H_1(C, \mathbf  Z)^-$. This gives a one to one
correspondence between $P(2^4)$ and the set of 
sub-lattices $H$ of $H_1(C, \bold Z)^{-}$ such that
\begin{enumerate}\item
$H=H^{\perp}$, 
\item
$H/(1-\rho)H_1(C, \bold Z)^{-}$ contains a vector $v$
with $q(v)\neq 0$.
\end{enumerate}
\end{proposition}

Under the correspondence of 
Proposition \ref{one to one cor  of index}, the lattice $L_1$
given in Proposition \ref{example good lattice}
corresponds to the partition
$\{\{1,2 \},\{3,4 \},\{5,6 \},\{7,8 \}\}$.

\subsection{Moduli spaces of branched coverings of $\mathbf  P^1$}
\label{moduli space of branced covering}
In this section, we give analytic descriptions of moduli
spaces of $4$-ple coverings of $\mathbf  P^1$ branching at
$8$ points.
Let $H_{std}$ be the module with a symplectic form and
the action of $\rho$ as in the last subsection.
An isomorphism $\phi : H_{std} \to H_1(C, \mathbf  Z)^-$
compatible with the intersection pairing and with the action of $\rho$ 
is called a marking of $C$. 
A pair $(C,\phi)$ 
is called a marked curve.
An isomorphism between two marked curves is an isomorphism 
between two curves which is compatible with the markings.
The set of isomorphic classes of marked curves is denoted by
$M_{marked}$. For an element $\sigma \in U(H_{std})$,
we define $\sigma (C,\phi)$ by $(C, \phi \circ \sigma)$.
By the commutative diagram

\setlength{\unitlength}{0.75mm}
\begin{picture}(200,45)(-40,-5)
\put(0,0){$H_{std}$}
\put(60,0){$H_1(C, \mathbf  Z)^-$}
\put(0,30){$H_{std}$}
\put(30,30){$H_{std}$}
\put(60,30){$H_1(C, \mathbf  Z)^-$}
\put(15,35){$\rho$}
\put(45,35){$\phi$}
\put(35,5){$\phi$}
\put(70,17){$\rho$}
\put(3,25){\line(0,-1){18}}
\put(4,25){\line(0,-1){18}}
\put(68,7){\vector(0,1){18}}
\put(12,30){\vector(1,0){16}}
\put(42,30){\vector(1,0){16}}
\put(10,4){\vector(1,0){45}}
\end{picture},

\noindent
we have $\rho (C, \sigma ) =(C, \sigma)$.

We define the complex ball $B(H_{std})$ associated to
$H_{std}$ as follows. The complex structure on the vector space 
$H_{std, \mathbf  R}$ is given by the action of $\mathbf  R(\rho)$
via the identification $\mathbf  R(\rho) \simeq \mathbf  C$
given by $\rho = -i$.
In the complex projective space $\mathbf  P(H_{std, \mathbf  R}^*)$,
the domain 
$$
B(H_{std})=\{h(v,*)\in \mathbf P(H_{std, \mathbf  R}^*) \mid h(v,v)<0 \}
$$
is isomorphic to a $5$-dimensional complex ball.
Then an element $\sigma \in U(H_{std})$ induces an
isomorphism of $B(H_{std})$.

For a marked curve $(C, \phi)$, we construct a point $p(C,\phi)$
in $B(H_{std})$ as follows. The complex structure
arising from $\rho$ is denoted by $\mathbf  R(\rho)$ 
to distinguish the usual complex structure.
By Proposition \ref{signature is 1.5},
the $(-i)$-eigenspace $H^0(C, \Omega^1)^{-,\rho=-i}$ of
$\rho$ in the space $H^0(C, \Omega^1)^-$ is one dimensional.
By the definition of a Hodge structure, we have an $\mathbf  R$-isomorphism
$$
H_1(C, \mathbf  R)^- \to Hom_{\mathbf  C}(H^0(C, \Omega^1)^-, \mathbf  C).
$$
By composing the map 
$$
Hom_{\mathbf  C}(H^0(C, \Omega^1)^-, \mathbf  C) \to
Hom_{\mathbf  C}(H^0(C, \Omega^1)^{-,\rho=-i}, \mathbf  C),
$$
and
$\displaystyle
H_{std,\mathbf  R} \overset\phi\simeq H_1(C, \mathbf  R)^-, 
$
we have an $\mathbf  R$-linear map
\begin{equation}
\label{period map as a linear form} 
H_{std,\mathbf  R} \to Hom(H^0(C, \Omega^1)^{-,\rho=-i}, \mathbf  C).
\end{equation}
It is easy to see that if we consider the complex structure
on $H_{std,\mathbf  R}$ by the action of $\mathbf  R(\rho)$,
the $\mathbf  R$-linear map (\ref{period map as a linear form}) 
is linear for the complex structures via the isomorphism $\mathbf  R(\rho)
\simeq \mathbf  C$ given by $\rho = -i$.
This linear form defines a point
in $B(H_{std})$. The corresponding point is denoted by $p(C, \phi)$.
This map 
$p: M_{marked} \to B(H_{std})$ is called a period map.
More explicitly, this map is given as follows. Let $\omega$ be
a basis of one dimensional complex space of
$H^0(C, \Omega^1)^{-,\rho=-i}$. 
Then the map (\ref{period map as a linear form})
$H_{std}\otimes \mathbf  R \to \mathbf  C$ 
is given by
$\gamma \mapsto\int_{\phi(\gamma)}\omega$.
By using the dual basis $A_j^*$ of $A_j$ over $\mathbf  R(\rho)$, this map
is expressed as $\sum_{j=1}^6(\int_{\phi(A_j)}\omega)A_j^*$.
This period map is holomorphic with respect to the parameters
$x_1, \dots ,x_8$.

\begin{theorem}[Terada, Deligne-Mostow]
The map $p$ is an open embedding and the complement of the
image is a proper analytic subset.
The natural actions of $U(H_{std})$ on $M_{marked}$ and
$B(H_{std})$ are compatible.
\end{theorem}

Let $M_{8pts}$ (resp. $M_{unord}$) be the set of 
the isomorphism classes of ordered
(resp. unordered)
set of distinct $8$ points in $\mathbf  P^1$.
The set $M_{8pts}$ (resp. $M_{unord}$) has a natural 
structure of an algebraic variety
which is isomorphic to $((\mathbf  P^1)^8- \Diag)/PGL(2, \mathbf  C)$
(resp. $((\mathbf  P^1)^8- \Diag)/\frak S_8 \times PGL(2, \mathbf  C)$),
where $\Diag = \{(x_j)|x_p =x_q\text{ for some }p<q\}$.
By corresponding a marked curve $(C,\phi) \in M_{marked}$ to 
the set of branching points $\{x_1, \dots ,x_8\}$ of 
$C\to \mathbf P^1$,  we have a morphism:
$$
M_{marked} \to M_{unord}.
$$

\begin{proposition}
Via the open immersion $p$, $M_{unord}$ is identified with an open set
of $B(H_{std})/U(H_{std})$. Moreover the covering
$M_{8pts}$ of $M_{unord}$ is identified with an open set
of $B(H_{std})/\Gamma(i+1)$,
where
$$
\Gamma (i+1) = \{ g \in U(H_{std}) \mid
 g \equiv 1\ {\rm mod }\ (1+\rho)H_{std} \}.
$$
\end{proposition}

We define complex reflections $M_{p,p+1} \in U(H_{std})$ for $p=1, \dots ,7$.
We choose
an initial point $X=(x_1, \dots ,x_8) \in M_{8pts}$ such that
$x_j \in \mathbf  R$ $(j=1,\dots,8)$ and $x_1 <\cdots < x_8$ (see figure in 
\S \ref{homology of covering br at 8 pts}).
The image of $X$ in $M_{unord}$ is denoted by $\overline{X}$.
We consider a path  
$M_{p,p+1}=(M_{p,p+1,j})_{j=1, \dots ,8} :[0,1] \to M_{8pts}$ by
\begin{align*}
& M_{p,p+1,j}(t)= x_j \text{ for }j\neq p, p+1, \\
& M_{p,p+1,p}(t)= \frac{x_{p+1}+x_p}{2}-\frac{x_{p+1}-x_{p}}{2}\mathbf  e(t/2),
 \\
& M_{p,p+1,p+1}(t)= \frac{x_{p+1}+x_p}{2}+\frac{x_{p+1}-x_{p}}{2}\mathbf  e(t/2),
\end{align*}
where $\mathbf e(t)=\exp(2\pi i t)$.
Then $M_{p,p+1}$ defines a closed path in $M_{unord}$ with
the base point $\overline{X}$.
By fixing the point $X$, we define the monodromy action $M_{p,p+1}$
of $H_{std}$ as follows. Let $C$ be the curve defined by the equation
(\ref{equation of branched curve}), where $x_1, \dots ,x_8$ are 
the coordinates of $X$. 
A basis $A_1, \dots, B_6$ of $H_1(C, \mathbf  Z)^-$
defined in \S \ref{homology of covering br at 8 pts}
gives a marking $H_{std}\to H_1(C, \mathbf  Z^-)$.
We consider the lifting $\tilde M_{p,p+1}$ 
of the path $M_{p,p+1}$ beginning from the point in $M_{marked}$
corresponding to the pair $(C, \phi)$. Then the end point 
$(C, \phi')$ of $\tilde M_{p,p+1}$ is a lifting 
$\overline X$ of $X$ in $M_{marked}$. The composite map
$\phi^{-1}\circ \phi'$ 
$$
H_{std} \overset{\phi'}\longrightarrow H_1(C, \mathbf  Z)^- 
\overset{\phi}\longleftarrow H_{std}
$$
is denoted by $M_{p, p+1}$. 
Since the paring and the action of $\rho$ are preserved
in the family $H_1(C_{M_{p,p+1}(t)}, \mathbf  Z)^-$,
($t \in [0,1]$),
we have $M_{p,p+1} \in U(H_{std})$.
Since $\Gamma (i+1)$ is a normal subgroup of $U(H_{std})$,
The covering 
$B(H_{std})/\Gamma (i+1) \to B(H_{std})/U(H_{std})$ is 
a Galois covering and
$$
Gal(B(H_{std})/\Gamma (i+1) \to B(H_{std})/U(H_{std}))
\simeq U(H_{std})/ \Gamma (i+1) \simeq O_6^+(2).
$$
By chasing the action of $M_{p,p+1}$ on $H_1(C, \mathbf  Z)^-$,
we have the following lemma.
\begin{lemma}
\label{O6+(2) to S8}
There exists an isomorphism $O_6^+(2) \simeq \frak S_{8}$
such that the image of $M_{p,p+1} \in U(H_{std})$ is
the transposition $(p,p+1)$ of $p$ and $p+1$.
Under the isomorphism
$$
B(H_{std})/\Gamma (i+1) \simeq M_{marked},
$$
the action of $\frak S_8 \subset Aut (M_{marked})$ is induced by
$\sigma^* (x_j)=x_{j\sigma}$.
Here the group $\frak S_8$ acts on the set $\{1, \dots ,8\}$ from
the right.
The action of $M_{p,p+1}$
is a complex reflection for the root $A_p$ with the eigenvalue $-\rho$, 
i.e., $M_{p,p+1}$ is characterized by 
$$
M_{p,p+1}(v)=\begin{cases}
-\rho ( A_p ) & \text{ if }v =A_p, \\
v & \text{ if }(v, A_p)=0. 
\end{cases}
$$
\end{lemma}
The isomorphism $O_6^+(2) \simeq \frak S_{8}$ in Lemma
\ref{O6+(2) to S8} induces a homomorphism $\pi:U(H_{std}) \to \frak S_8$. 

We define an inclusion $B(H_{std})$ to the Siegel upper half space 
$\frak H_6$ of degree $6$ by using a good basis
$(L,\Sigma_L)$ of $H_{std}$ as follows.
By the Poincare duality, 
we have $h(y_1, y_2)=0$ for elements
$y_1 \in \ker(H_{std}, Hom(H^0(C, \Omega^1)^{-,\rho=-i}, \mathbf  C))$
and
$y_2 \in \ker(H_{std}, Hom(H^0(C, \Omega^1)^{-,\rho=i}, \mathbf  C))$, 
where $H^0(C, \Omega^1)^{-,\rho=i}$ is the 
$i$-eigenspace of $\rho$ in the space $H^0(C, \Omega^1)^-$. 
Let $H$ be the corresponding Hodge structure of $p \in B(H_{std})$.
Let $\Sigma_{L}=\{ a_1, \dots ,a_6, b_1, \dots ,b_6\}$ and
$$
H \otimes \mathbf  C \simeq H^{(1,0)} \oplus H^{(0,1)}
$$ 
be the Hodge decomposition of $H$. We choose a basis  
$\omega_1, \dots ,\omega_g$ of $H^{(1,0)}$ such that 
$\int_{b_j}\omega_k= \delta_{jk}$. We put 
$\tau_{jk}=\int_{a_j}\omega_k$.
Then by the definition of
polarized Hodge structure, $\tau =(\tau_{jk})_{jk}$ is an
element of $\frak H_6$.

\subsection{Level $2$ structure and exponent $2$ covering of
configuration space}
\label{level2 str and 2 expo. cov}

Let $\mathbf  C(x_1,\dots ,x_8)$ be the rational function field
of $x_1, \dots, x_8$ over $\mathbf  C$. On this field, 
the groups $PGL(2, \mathbf  C)$ and $\frak S_8$ act by
$$
g(x_j) = \frac{a x_j +b}{c x_j +d},\quad 
\sigma (x_j)=x_{j\sigma}
$$
for $g\in PGL(2, \mathbf  C)$ and $\sigma \in \frak S_8$, respectively. 
As in the last section the group $\frak S_8$ acts on the 
set $\{1, \dots ,8\}$ from the right.
(Note that the action of $\frak S_8$ on the space $M_{marked}$
is covariant.)
Let $K$ be the fixed subfield of 
$M=\mathbf  C(x_1,\dots ,x_8)$ under the action of $PGL(2, \mathbf  C)$.
Since the action of $\frak S_8$  commutes with that of $PGL(2, \mathbf  C)$, 
$\frak S_8$ acts on $K$.
Let $L$ be the fixed subfield of $K$ under the action of $\frak S_8$.
Then $K$ and $L$ are equal to the function fields of $M_{8pts}$
and $M_{unord}$, respectively.
The field $K$ is generated by the cross-ratios 
$$ 
\lambda_j=\frac{(x_3-x_1)(x_j-x_2)}{(x_2-x_1)(x_j-x_3)}$$
of
$\{x_1,x_2,x_3,x_j \}$ for $j=4, \dots ,8$ over $\mathbf  C$.

Let $\tilde K$ be the
algebraic closure of $K$ in 
$\tilde M=\mathbf  C(x_j,\sqrt{x_j-x_k})_{j \neq k}$.
Since the extension $\tilde M/M$ is a Galois extension,
so is $\tilde K/ K$. The extensions $M$ and
$\tilde K$ are linearly independent over K, therefore the restriction
map 
$$
Gal (\tilde M/ M) \to Gal (\tilde K/K)
$$
is surjective, and $\tilde K$ is generated by
$f=\prod_{j <k}\sqrt{x_j- x_k}^{a_{jk}}$
such that $f^2$ is an element of $K$. Thus
$\tilde K$ is generated by $\sqrt{\lambda_j}$, 
$\sqrt{\lambda_j-1}$ for $j=4, \dots , 8$, and 
$\sqrt{\lambda_j-\lambda_k}$ for $4 \leq j < k \leq 8$
and $\tilde K$ is a Galois extension
of $L$.  The inclusions of fields
$L \subset K \subset \tilde K$ imply the following exact sequence
of groups:
$$
1 \to N \to Gal(\tilde K/L) \to \frak S_8 \to 1,
$$
where $N =Gal (\tilde K/K) \simeq (\mathbf  Z/2\mathbf  Z)^{20}$.

We compare this Galois extension with the analytic description of 
the corresponding moduli space of $4$-ple coverings of
$\mathbf  P^1$ branching at eight points.

\begin{proposition}
Let $\tilde M_{8pts}$ be the normalization of the 
$M_{8pts}$ in $\tilde K$. This variety $\tilde M_{8pts}$
is identified with an open subset of $B(H_{std})/\Gamma (2)$
via the period map $p$,
where
$$
\Gamma (2) = \{ g \in U(H_{std}) \mid
 g \equiv 1\ {\rm mod }\ 2H_{std} \}.
$$
Via this isomorphism, we have
$$Gal(\tilde K/L) \simeq U(H_{std})/\Gamma(2)\cdot \langle i \rangle,$$
where $\langle i \rangle$ is the cyclic group generated by $i.$

\end{proposition}
\begin{proof}

We already know that
$$
U(H_{std}) / \Gamma (i+1)\cdot \langle i \rangle \simeq \frak S_8.
$$
The group
$\Gamma (i+1)\cdot \langle i \rangle /\Gamma (2)\cdot \langle i \rangle 
\simeq (\bold Z/2\bold Z)^{20}$
is generated by $M_{jk}^2$. By restricting the action of
$U(H_{std})$ to $\tilde K$, we have a map
\begin{equation}
\label{partial isom}
\Gamma (i+1)\cdot \langle i \rangle /\Gamma (2)\cdot \langle i \rangle 
 \to
N.
\end{equation}
We study the action of $M_{jk}$ on 
$\{\sqrt{x_p -x_q}\}$.
We assign each algebraic
function $\sqrt{x_j-x_k}$ on $(\bold P^1)^8 - \Diag$
a branch as follows. 
Let $X$ be the initial point in $M_{8pts}$ as in the last section.
For $j<k$, $\sqrt{x_k - x_j}$
denotes the branch of the algebraic function on 
$(\bold P^1)^8 - \Diag$
so that it takes a positive real value at $X$. 
The analytic continuation of the function $\sqrt{x_k - x_j}$ 
along the path $M_{p,p+1}$ is given as
$$
M_{p,p+1}(\sqrt{x_k -x_j})=
\begin{cases}
i \sqrt{x_k -x_j}
& \text{ if $j=p$ and $p+1 =k$}, \\
\sqrt{x_{k \sigma} -x_{j \sigma }} 
& \text{ otherwise }, \\
\end{cases}
$$
where $\sigma$ is the transposition $(p,p+1)$.
As a consequence, we have 
$$
M_{p,p+1}^2(\sqrt{x_k -x_j})=
\begin{cases}
- \sqrt{x_k -x_j}
& \text{ if $j=p$ and $p+1 =k$}, \\
\sqrt{x_{k  } -x_{j  }} 
& \text{ otherwise }. \\
\end{cases}
$$
By looking at the action of $M_{jk}^2$ on the set
$\{ \sqrt {\lambda_p}, \sqrt{1 -\lambda_p}, \sqrt{\lambda_p - \lambda_q}
\}$, the homomorphism (\ref{partial isom})
is an isomorphism.

\end{proof}

Let $g \in U(H_{std})$ and 
$r_1 =\{\{1,2\},\{3,4\},\{5,6\},\{7,8\}\}\in P(2^4)$.
We put
$r_1\pi(g)=\{\{j_1,j_2\},\{j_3,j_4\},\{j_5,j_6\},\{j_7,j_8\}\} \in P(2^4)$,
where $\pi:U(H_{std}) \to \frak S_8$. 
We assume that $j_p < j_{p+1}$ for $p=1,3,5,7$. Then there exists 
a $4$-th root of the unity $arg(g)$, called the argument of $g$, such that
\begin{align*}
& g(\sqrt{(x_2-x_1)(x_4-x_3)(x_6-x_5)(x_8-x_7)}) \\
= &
arg(g)\sqrt{(x_{j_2}-x_{j_1})(x_{j_4}-x_{j_3})
(x_{j_6}-x_{j_5})(x_{j_8}-x_{j_7})}.
\end{align*}

\section{Theta function of standard principal sub-Hodge structure $L_1$}
\label{Theta principal}
\subsection{Abel-Jacobi map and the order of zero}
\label{Abel-Jacobi map and order of zero}
Let $p_1, \dots , p_8$ be the ramification points
of the smooth curve $C$ defined by (\ref{curve C})
above $x_1, \dots ,x_8$, respectively.
Let $jac$ be the Abel-Jacobi map 
\begin{eqnarray*}
jac:C &\longrightarrow  &J(C) =Hom(H^0(C, \Omega^1), \mathbf  C)/
H_1(C, \mathbf  Z)
 \\ 
p &\mapsto  & 
\text{ the linear function } \int_{p_1}^p \text{ on }H^0(C, \Omega^1)
\text{ defined by }\\
& & \int_{p_1}^p:  H^0(C, \Omega^1)\ni \omega\mapsto 
\int_{p_1}^p  \omega\in \mathbf C. 
\end{eqnarray*}
The endomorphism $(1-\rho^2):J(C) \to J(C)$ of $J(C)$
factors through the natural inclusion $\kappa : Prym(C) \to J(C)$,
i.e. there is a morphism  $\alpha :J(C) \to Prym(C)$ such
that $\kappa \circ \alpha =1-\rho^2$. Let $J(C)^+$ be the
the connected component of the kernel of $\alpha$. Then it is easy
to see that $J(C)/J(C)^+$ is isomorphic to $Prym(C)$
and that the morphism 
$$
J(C)/J(C)^+ \simeq Prym(C) \to Prym(C)
$$
induced by the morphism $\alpha$
corresponds to the index finite group $(1-\rho^2)H_1(C, \mathbf  Z)=
(1-\rho)H_1(C, \mathbf  Z)^-$ of $H_1(C, \mathbf  Z)^-$.
As a consequence, we have the following sequence of
morphisms:
$$
J(C) \to J(C)/J(C)^+ \simeq Prym(C) \to Prym(C) \to J(C). 
$$
The composite map 
$$
C \overset{jac}\rightarrow J(C) \rightarrow J(C)/J(C)^+ \simeq Prym(C)
$$ 
is denoted by $jac^-$.
Let $L_1$ be the principal sub-Hodge structure defined 
in \S\ref{polarization of prym var}.
Then $A_{L_1} = \mathbf  C^6/L_1$ is a principally
polarized abelian variety.
The inclusions 
$$
(1-\rho)H_1(C, \mathbf  Z)^- \subset L_1 \subset H_1(C, \mathbf  Z)^-
$$
induce homomorphisms of abelian varieties
$$
Prym (C) \overset{\pi_1}\longrightarrow A_{L_1} \overset{\pi_2}\longrightarrow
Prym (C).
$$
We define the theta function $\vartheta_m(\Sigma_1,z)$ for 
the good basis $\Sigma_1$ defined in \S\ref{polarization of prym var} 
with the characteristic   
$m=(m', m'') \in \mathbf  Q^{12}$ by 
$$
\vartheta_m(\Sigma_1,z) =
\sum_{\xi \in \mathbf  Z^6}\mathbf  e(\frac{1}{2}(\xi +m')\tau\;^t(\xi+m')+
(z+m'')\;^t(\xi+m')),
$$
where $z \in \mathbf  C^6$ and
$\tau = (\tau_{ij})_{ij}\in \frak  H_6$ is the normalized period matrix 
for the good basis $\Sigma_1$ defined in the last section.

Let $\iota :C \to A_{L_1}$ be the composition $\pi_1 \circ jac^-$.
By Proposition \ref{example good lattice} and
\ref{one to one cor  of index}, 
we have the following lemma.
\begin{lemma}
\begin{enumerate}
\item
The image of each ramification point $p_1, \dots ,p_8$
under the map $\iota$ is contained in the set of 
$(1-\rho)$-torsion points of $Prym(C)$.
\item
The image of $\pi_1$ is identified with the quotient of
$Prym(C)$ by the group generated by 
$ jac^{-}(p_2),jac^{-} (p_3)-jac^{-} (p_4),jac^{-} (p_5)-jac^{-} (p_6)$.
\item
Using the good base $\Sigma_{L_1}$ of $L_1$, we have
$$
\iota(p_{k})\equiv\frac{1}{2}(\xi_k, \xi_kU) \text{ mod } L_1,
$$
where
$\xi_1=\xi_2=0$,
$\xi_3 =\xi_4=(1,1,0,0,0,0)$, $\xi_5=\xi_6=(1,1,1,1,0,0)$
and $\xi_7=\xi_8=(1,1,1,1,1,1)$.
\end{enumerate}
\end{lemma}

We study the order of zero of 
the pull back of $\vartheta_m(\Sigma_1, z)$ by $\iota$.
Let $\tilde C$ be the universal covering of $C$ and
we choose a base point $\tilde q_1$ of $\tilde C$ as a lifting of $p_1$. 
Then we have a lifting 
$\tilde\iota: \tilde C \to (H^0(C, \Omega^1)^-)^*$ 
of $\iota$ by sending $\tilde p_1$ to $0$.
Let $\omega_1, \dots \omega_6$ be the normalized basis of $H^0(C, \Omega^1)^-$ 
with respect to $\Sigma_1$. Via the isomorphism
$$(H^0(C, \Omega^1)^-)^*\ni \gamma \mapsto  
(\int_{\gamma}\omega_1,\dots,\int_{\gamma}\omega_6)\in \mathbf  C^6,$$ 
$\tilde \iota$ is identified with $\tilde C \to \mathbf  C^6$.
We define a map $F_{m} :\tilde C \to \mathbf  C$ by 
$$
\tilde C \ni \tilde p \mapsto 
F_m(\tilde p)=\vartheta_{m}(\Sigma_1, \tilde \iota (\tilde p))
\in C. 
$$
Since $\vartheta_m(\tau,z)$ is a non-zero section of a line bundle
$\Cal L_m$,
the order of zero at the lifting $\tilde p$ of $p$
depends only on the point $p$ in $C$.
It is called the order of zero at $p$ and denoted by $ord_p(F_m)$.
\begin{proposition}
The total sum of the order $ord_p(F_m)$ of $p$
on $C$ is 12.
\end{proposition}
\begin{proof}
Use a similar argument in \cite{MT} just after Proposition 4.9.
\end{proof}
The next proposition is fundamental for determining the
distribution of zeros of $F_m$.
\begin{proposition}
\label{order of theta mod 4}
Let $m=\frac{1}{2}(\mu, \mu U)$ ($\mu \in \bold Z^6$) be an element of 
$\frac{1}{1-\rho}L_1$. ( See Remark \ref{torsion}).
Let $\xi_j$ be an element of $\bold Z^6$ such that
$\iota(\tilde p_j)\equiv\frac{1}{2}(\xi_j, \xi_j U)$ (mod $\mathbf  Z^{12}$),
where $\tilde p_j$
is a lifting of $p_j$ to $\tilde C$.
We put $q=\mu + \xi_j$.
Then $ord_{p_j}(F_m)$ is equal to
$-qU\;^tq$ modulo $4$.
\end{proposition}
\begin{proof}
Let $z$ be the coordinate for the universal covering
of $A_{L_1}$.
Since the point $p_j$ is fixed under the action of $\rho$,
$\rho z = z + l$ ($l \in L_1$).
By the transformation formula in p.85,
\cite{I}, we have $F(\rho z) = u(z) F(z)$, where
$$
\lim_{z \to \tilde p_j} u(z) = \bold e( \frac{-qUq}{4}).
$$
Therefore the order of $F_m$ at $p_j$ is congruent to
$-qUq$ mod $4$.
\end{proof}
Let $m_k=\frac{1}{2}(\mu_k,\mu_kU)$ ($k=1,\cdots ,4$), where
\begin{align}
\label{character list}
&\mu_1=(0,0,0,0,0,0),\quad\mu_2=(0,0,1,1,1,1), \\
\nonumber
&\mu_3=(1,1,0,0,1,1),\quad\mu_4=(1,1,1,1,0,0).
\end{align}
Then by Proposition \ref{order of theta mod 4}, the
table of
$ord_{p_j}(F_{m_k})$ (mod $4$) is given by
\begin{equation}
\label{order of zeros}
\begin{matrix}
 p_1 & p_2 & p_3 & p_4 & p_5 & p_6 & p_7 & p_8 \\
 0 & 0 & 2 & 2 & 0 & 0 & 2 & 2  \\
\end{matrix}
\end{equation}
for $k=1, \dots ,4$.

\subsection{Determination of extra zeros of theta functions}
\label{determination}

Let $\mu_1,\dots, \mu_4$ and $m_1, \dots  ,m_4$ be as in \S 
\ref{Abel-Jacobi map and order of zero}.
Then the sum of the known zeros of $F_{m_k}(\tilde p)$
is $8$ for $k=1, \dots ,4$.
Since the action of $\rho$ on the curve $\tilde C$ and 
$H^0(C,\Omega^1)$ is compatible, the remaining $12-8=4$ zeros are
stable under the action of $\rho$ and if its support contains
one of $p_j$, then its multiplicity should be $4$ 
by the modulo $4$ condition.

\begin{proposition}
\label{invariance of diff char} 
The function $R_{jk}=F_{m_j}(\tilde p)/F_{m_k}(\tilde p)$
is a rational function of $p=(z,w)\in C$ for $1 \leq j,k \leq 4$.
Moreover 
$R_{jk}(\tilde p)=\displaystyle c\cdot \frac{z-s}{z-t}$,
with some constants $s,t$ and $c \neq 0$. 
\end{proposition}
\begin{proof}
Since the image of the fundamental group of $C$ in
$L_1$ is equal to $(1-\rho)H_1(C, \bold Z)^-$, if
the theta functions $\vartheta_{m_j}(\tau,z)$ and
$\vartheta_{m_k}(\tau, z)$
have the same quasi-periodicity, the quotient $F_{m_j}/F_{m_k}$
is a rational function on $C$. By comparing the zeros of the
numerator and the denominator of $R_{jk}$
(see table (\ref{order of zeros})), 
we have the proposition.
\end{proof}
\begin{proposition}
\label{st equation}
\begin{enumerate}
\item
In the expression of the rational function 
$R_{13}=c \cdot\displaystyle\frac{z-s}{z-t}$,
$s$ and $t$ are determined by the
equation:
\begin{equation}
\label{equation of t and s}
\frac{x_1-s}{x_1-t} + \frac{x_2-s}{x_2-t} =0,\quad 
\frac{x_5-s}{x_5-t} + \frac{x_6-s}{x_6-t} =0.
\end{equation}
\item
The rational function $R_{12}$ on $C$ is a constant.
\end{enumerate}
\end{proposition}
\begin{proof}
By Proposition \ref{invariance of diff char},
we have
$$
c\cdot\frac{x_j-s}{x_j-t} =
R_{13}(\tilde p_j) =
\frac{F_{m_1}(\tilde p_j)}{F_{m_3}(\tilde p_j)}.
$$
On the other hand, we have
$$
\frac{F_{m_1}(\tilde p_1)}{F_{m_3}(\tilde p_1)} = 
-\frac{F_{m_1}(\tilde p_2)}{F_{m_3}(\tilde p_2)}
$$
by the quasi-periodicity of theta functions, thus we have the
statement (1).
We can prove the statement (2) similarly.
\end{proof}

\begin{proposition}
Let $\vartheta_k=\vartheta_{m_k}(\tau)$. Then we have
$$
\frac{(\vartheta_{2}+i\vartheta_{3})^2(\vartheta_{1}-i\vartheta_{4})^2}
{4\vartheta_{1}^2\vartheta_{3}^2} =
\frac{(x_1-x_5)(x_2-x_6)}{(x_1-x_2)(x_5-x_6)}.
$$
\end{proposition}
\begin{proof}
By the definition of theta constants and $R_{jk}(\tilde p)$, we have
\begin{align*}
R_{13}(p_1)=\frac{\vartheta_{1}}{\vartheta_{3}}=
c\cdot\frac{x_1-s}{x_1-t}, \quad
R_{13}(p_5)=-\frac{\vartheta_{4}}{\vartheta_{2}}=
c\cdot\frac{x_5-s}{x_5-t}.
\end{align*}
The equality $R_{12}(p_1)=R_{12}(p_5)$
implies
\begin{align*}
\frac{\vartheta_{1}}{\vartheta_{2}}=
\frac{\vartheta_{4}}{\vartheta_{3}}.
\end{align*}
By computing $R_{13}(p_1)/R_{13}(p_5)$, we have
\begin{align}
\label{ratio of special value of rat. fn.}
\frac{(x_1-s)(x_5-t)}{(x_1-t)(x_5-s)} = 
-\frac{\vartheta_{1}\vartheta_{2}}
{\vartheta_{3}\vartheta_{4}} 
= 
-\frac{\vartheta_{1}^2}
{\vartheta_{4}^2} =
-\frac{\vartheta_{2}^2}
{\vartheta_{3}^2} 
\end{align}
and
\begin{align*}
(1+\frac{\vartheta_{4}^2}{\vartheta_{1}^2})
(1+\frac{\vartheta_{1}^2}{\vartheta_{4}^2})
= & \frac{(\vartheta_{1}^2+\vartheta_{4}^2)^2}
{\vartheta_{4}^2\vartheta_{1}^2} =
\frac{\vartheta_{2}^2+\vartheta_{3}^2}
{\vartheta_{2}\vartheta_{3}}\cdot
\frac{\vartheta_{1}^2+\vartheta_{4}^2}
{\vartheta_{4}\vartheta_{1}} \\
= &\frac{(\vartheta_{2}+i\vartheta_{3})(\vartheta_{2}-i\vartheta_{3})
(\vartheta_{1}+i\vartheta_{4})(\vartheta_{1}-i\vartheta_{4})}
{\vartheta_{1}\vartheta_{2}\vartheta_{3}\vartheta_{4}} \\
= &
\frac{(\vartheta_{2}+i\vartheta_{3})^2(\vartheta_{1}-i\vartheta_{4})^2}
{\vartheta_{1}^2\vartheta_{3}^2}.
\end{align*}
Here we used the relation 
$
(\vartheta_{2}-i\vartheta_{3})(\vartheta_{1}+i\vartheta_{4})=
(\vartheta_{2}+i\vartheta_{3})(\vartheta_{1}-i\vartheta_{4})$.
On the other hand, by
(\ref{equation of t and s}) and
(\ref{ratio of special value of rat. fn.}), we have
$$
\frac{1}{4}(1+\frac{\vartheta_{4}^2}{\vartheta_{1}^2})
(1+\frac{\vartheta_{1}^2}{\vartheta_{4}^2})
=
\frac{(x_1-x_5)(x_2-x_6)}{(x_1-x_2)(x_5-x_6)}.
$$
\end{proof}

\subsection{An application of quadratic theta relation}
\label{application of quadr. rel of theta}
In this section, we fix a principally polarized sub-Hodge structure
$L=L_1$ and study quadratic relations between theta constants.
We recall the quadratic relation between theta functions in \cite{I}.
For the next proposition, see p.139, \cite{I}. 
\begin{proposition}
Put 
$$
n_1 =(n_1',n_1'')= \frac{1}{2}(m_1 + m_2), \quad 
n_2 =(n_2',n_2'')=\frac{1}{2}(m_1-m_2),
$$
for $m_1=(m_1', m_1''),\ m_2=(m_2',m_2'')\in \mathbf  Q^{12}$. 
Let $S$ be a complete set of representatives of
$(\frac{1}{2}\mathbf  Z)^6/\mathbf  Z^6$. We have
$$
\vartheta_{m_1}(\tau)\vartheta_{m_2}(\tau)=
\frac{1}{2^6}\sum_{a'' \in S}\mathbf  e(-2m_1^{' t}a'')
\vartheta_{2n_1',n_1''+a''}(\frac{\tau}{2})
\vartheta_{2n_2',n_2''+a''}(\frac{\tau}{2}).
$$
\end{proposition}

We apply this formula to 
$$
n_1=(\frac{1}{2}v_1-\frac{1}{4}U_0, \frac{1}{2}U_0) 
\quad n_2=(-\frac{1}{4}U_0,0), 
$$
where $v_1\in \mathbf  Z^6$, and
we use Notation \ref{beginig notation}
for $U_0=(1,1,0,0,1,1)$.  
Replace $\tau$ by $\tau+U$, then we have
\begin{align}
\label{quadr relation}
& \vartheta_{\frac{1}{2}v_1-\frac{1}{2}U_0,\frac{1}{2}U_0}(\tau+U)
\vartheta_{\frac{1}{2}v_1,\frac{1}{2}U_0}(\tau+U) \\ 
= & \frac{1}{2^6}\sum_{a'' \in S}\mathbf  e(-(v_1-U_0)^{t}a'')
\vartheta_{v_1-\frac{1}{2}U_0,\frac{1}{2}U_0+a''}(\frac{1}{2}(\tau+U))
\vartheta_{-\frac{1}{2}U_0,a''}(\frac{1}{2}(\tau+U)).
\nonumber
\end{align}

We assume that $\tau$ is the normalized period matrix of the
principally polarized Hodge structure $L$ with respect to
the symplectic basis $\Sigma_1$. Then we have 
\begin{equation}
\label{rho action and period matrix}
(\tau U)^2 = -I.
\end{equation}
By applying the transformation formula 
$$
\vartheta_{m',m''}(\tau +U)=
\mathbf  e(-\frac{1}{2}m' U\;^t m' + \frac{1}{2}m^{'t}U_0)
\vartheta_{m',m''+m'U+\frac{1}{2}U_0}(\tau)
$$
to the left hand side of (\ref{quadr relation}), we have
\begin{align}
\label{left hand side of q.rel}
& \vartheta_{\frac{1}{2}v_1-\frac{1}{2}U_0,\frac{1}{2}U_0}(\tau+U)
\vartheta_{\frac{1}{2}v_1,\frac{1}{2}U_0}(\tau+U)  \\
= &
\mathbf  e(-\frac{1}{4}v_1U\;^tv_1+\frac{3}{4}v_1\;^tU_0-\frac{3}{8}U_0\;^tU_0)
\vartheta_{\frac{1}{2}(v_1-U_0),\frac{1}{2}(v_1-U_0)U}(\tau)
\vartheta_{\frac{1}{2}v_1,\frac{1}{2}v_1U}(\tau).
\nonumber 
\end{align}

To compute the right hand side, we apply the transformation 
formula in p.85, \cite{I},
\begin{equation}
\label{igusa transformation formula}
\vartheta_{m^\#}(\tau^{\#})=
\det (C \tau +D)^{1/2}\cdot u \cdot \vartheta_m(\tau),
\end{equation}
for
$$
\sigma =\left(\begin{matrix}
A & B \\ C & D \end{matrix}\right) =
\left(\begin{matrix}
0 & U \\ -U & I \end{matrix}\right) 
$$
and $m=(-(a''+\frac{1}{2}U_0)U,v_1U+a''+\frac{1}{2}U_0)$
(resp. $m=(-a''U,a'')$).
Here 
$\tau^\#=(A\tau +B)(C\tau +D)^{-1}$ is equal to $\frac{1}{2}(\tau +U)$
and $\det (C\tau +D)=-8$
by the relation (\ref{rho action and period matrix}).
The theta characteristic
$$m^\#=m\cdot \sigma^{-1}+ \frac{1}{2}((C\;^tD)_0, (A\;^tB)_0)$$
is equal to $(v_1-\frac{1}{2}U_0, \frac{1}{2}U_0+a'')$
(resp. $(-\frac{1}{2}U_0, a'')$). 
We fix a branch of $\det (C\tau +D)^{1/2} = \sqrt 8 i$ once and for all.
We compute the constant $u$ in the formula
(\ref{igusa transformation formula})
which depends only on $m$.

\begin{definition}
We define a non-zero complex number $c(a,b)$ by
\begin{equation}
\label{def of c(a,b)}
\vartheta_{a,b}(U(-U\tau+I)^{-1}, z^\#)=
c(a,b)(-8)^{1/2}\vartheta_{c,d+\frac{1}{2}U_0}(\tau, z),
\end{equation}
where $z^\#=z(-U\tau+I)^{-1}$ and $c=-bU, d=aU+b$.
\end{definition}
\begin{proposition}
\label{root of unity in right hand side}
\begin{enumerate}
\item
$c(a,b)/c(0,0)=\mathbf  e(\frac{1}{2}bU\;^tb+a\;^tb +\frac{1}{2}b\;^tU_0)$.
\item
$c(0,0)^2=1$.
\end{enumerate}
\end{proposition}
\begin{proof}
1. This is the direct consequence of the formula in p.85, \cite{I}.

2. Since $c(0,0)$ is independent of $\tau$, we evaluate the both sides of
(\ref{def of c(a,b)}) at $\tau = iI +U$. Then we have
\begin{align*}
\vartheta_0(iI)= &\det (-U\tau +I)^{1/2}c(0,0)
\vartheta_{0,\frac{1}{2}U_0}(iI +U) \\
= &(-1) \cdot \det (-i U)^{1/2}c(0,0)
\vartheta_{0,0}(iI ).
\end{align*}
Since $\vartheta_0(iI) \neq 0$, we have  $c(0,0)^2=1$.
\end{proof}

We can compute the right hand side of (\ref{quadr relation}) by 
Proposition \ref{root of unity in right hand side}. 
As a consequence we have the following theorem.
For an element $m \in \frac{1}{2} \mathbf  Z^{6}$, 
a representative of the class
of $m$ in $\frac{1}{2}\mathbf  Z^{6}/\mathbf  Z^{6}$ in 
$\{ 0, \frac{1}{2}\}^{6}$ is denoted by $<m>$.
\begin{theorem}
\label{comapre inside same L}
\begin{enumerate}
\item
Let $v_1$ be an element of $\bold Z^6$. Then we have
\begin{align*}
& 8 \mathbf  e(-\frac{1}{4}v_1U\;^tv_1 + \frac{3}{4}v_1\;^tU_0
-\frac{3}{8}U_0\;^tU_0)\cdot 
\vartheta_{\frac{1}{2}(v_1-U_0), \frac{1}{2}(v_1-U_0)U}(\tau)
\vartheta_{\frac{1}{2}v_1, \frac{1}{2}v_1 U}(\tau) \\
= & 
\sum_{a''}
\mathbf  e(a''U\;^ta''+\frac{3}{2}a^{''t}U_0-v_1\;^ta'') 
\vartheta_{a''U+\frac{1}{2}U_0, a''+\frac{1}{2}U_0}(\tau)
\vartheta_{a''U, a''}(\tau),
\end{align*}
where $a''$ runs over the complete set $S$ of representatives
in $\frac{1}{2}\bold Z^6/\bold Z^6$.
\item
Let $v_1 \in \{0,1\}^6$ and choose $b''$ as
$$
b''\equiv \frac{1}{2}v_1-\frac{1}{2}U_0,\quad 
b'' \in \{ 0, \frac{1}{2}\}^6.
$$
Then we have
\begin{align*}
& 8 \mathbf  e(-\frac{1}{4}v_1U\;^tv_1 + \frac{3}{4}v_1\;^tU_0
-\frac{3}{8}U_0\;^tU_0)\cdot 
\mathbf  e(\frac{1}{2}v_1\;^tU_0)
\vartheta_{b'', b''U}(\tau)
\vartheta_{\frac{1}{2}v_1, \frac{1}{2}v_1 U}(\tau) \\
= & 
\sum_{a''}
\mathbf  e(\frac{1}{2}a^{''t}U_0+v_1\;^ta'') 
\vartheta_{<a''+\frac{1}{2}U_0>U, <a''+\frac{1}{2}U_0>}(\tau)
\vartheta_{a''U, a''}(\tau).
\end{align*}
Moreover if $\frac{1}{4}v_1U\;^tv_1 \in \mathbf  Z$, then we have
\begin{align*}
&- 8 \mathbf  e( \frac{1}{4}v_1\;^tU_0)\cdot  
\vartheta_{b'', b''U}(\tau)
\vartheta_{\frac{1}{2}v_1, \frac{1}{2}v_1 U}(\tau) \\
= & 
\sum_{a''}
\mathbf  e(\frac{1}{2}a^{''t}U_0+v_1\;^ta'') 
\vartheta_{<a''+\frac{1}{2}U_0>U, <a''+\frac{1}{2}U_0>}(\tau)
\vartheta_{a''U, a''}(\tau).
\end{align*}
\end{enumerate}
\end{theorem}
We have the following corollary of Theorem \ref{comapre inside same L}.
\begin{corollary}
\label{cor to quadratic rel final}
Let $v_1 \in \{0,1\}^6$ such that
$v_1U\;^tv_1 \in 4\mathbf Z$, $v_1 \neq 0$, and
let $b''$ be the element in $\{0, \frac{1}{2}\}^6$ defined in 
Theorem \ref{comapre inside same L}.
Then we have
$$
\mathbf  e( \frac{1}{4}v_1\;^tU_0)\cdot  
\vartheta_{b'', b''U}(\tau)
\vartheta_{\frac{1}{2}v_1, \frac{1}{2}v_1 U}(\tau) 
+
\vartheta_{0, 0}(\tau)
\vartheta_{\frac{1}{2}U_0, \frac{1}{2}U_0}(\tau)=0. 
$$
\end{corollary}

\section{Comparison for theta constants of different lattices} 
\subsection{Translation vector arising from 
changing lattices}

Let $L$ be a principally polarized sub-Hodge structure
of $H_1(C, \mathbf  Z)^-$ and 
$$
\Sigma=\Sigma_L=
\{a_1, \dots, a_6, b_1, \dots ,b_6\}
$$ 
be a
good symplectic basis of $L$. 
Let $\{\alpha'_1, \dots, \alpha'_6, \beta'_1, \dots ,\beta'_6\}$
be the basis of $B$ defined in 
Proposition \ref{base for basic principal lattice}.
Put
\begin{align*}
&\alpha_j = \alpha'_j\ (\text{ for }j=1, \dots ,6), \\
&\beta_j=2\beta'_j\ (\text{ for }j=1, 2 ,3),\quad
\beta_j=\beta'_j\ (\text{ for }j=4,5,6),
\end{align*}
and $\Sigma_B=\{\alpha_1, \dots, \alpha_6, \beta_1, \dots ,\beta_6\}$.
Then the lattice $B$ generated by $\Sigma_B$ admits a principally polarized
sub-Hodge structure of $H_1(C, \mathbf  Z)^-$.
In this section we compare theta functions of $(B, \Sigma_B)$ and those of 
$(L, \Sigma_L)$.
The elements $a_j, b_j$ in $\Sigma_L$ are linear combinations
of $\Sigma_B$ as
\begin{align*}
& a_j=\sum_{j=1}^6a_{jk}\alpha_k+\sum_{j=1}^6b_{jk}\beta_k, \\
& b_j=\sum_{j=1}^6c_{jk}\alpha_k+\sum_{j=1}^6d_{jk}\beta_k,
\end{align*}
respectively. The column vector consisting of $\alpha_j$ (resp. $\beta_j$,
$a_j$ and $b_j$) for $j=1\dots6$ is denoted by $\alpha$ 
(resp. $\beta$, $\mathbf  a$ and $\mathbf  b$). 
We put $A=(a_{jk})$, $B=(b_{jk})$, $C=(c_{jk})$, $D=(d_{jk})$.
Then we see that 
$$
\sigma= 
\left(\begin{matrix}
A & B \\ C & D
\end{matrix}\right) 
\in Sp(6,\mathbf  Q).
$$
Let $p=(p_1, \dots ,p_6),q=(q_1, \dots ,q_6)$ be elements in $\mathbf  Q^6$.
We define two vectors $r=(r_1, \dots ,r_6)$ and $s=(s_1, \dots ,s_6)$
in $\mathbf  Q^6$ by
$$
(r,s)=(p,q)
\left(\begin{matrix}
A & B \\ C & D
\end{matrix}\right) .
$$
We have $(r,s)\;^t(\alpha,\beta)=(p,q)\;^t(\mathbf  a, \mathbf  b)$.
Then $(r,s)\;^t(\alpha,\beta)$ is an element in $H_1(C, \mathbf  Z)^-$ 
(resp. $(1-\rho)H_1(C, \mathbf  Z)^-$)
if and only if
$r\in \mathbf  Z^6$ and $s\in \frac{1}{2}\mathbf  Z^3 \oplus \mathbf  Z^3$
(resp.
$r\in 2\mathbf  Z^3 \oplus \mathbf  Z^3$ and $s\in \mathbf  Z^6$). 

\begin{proposition}
Let $e={\rm diag} (2,2,2,1,1,1)$. Then all entries of
$e\;^tB$,$e\;^tD$ are integers.
\end{proposition}

Let $\tau$ and $\tau^\#$ be the normalized period matrix of $B$
and $L$ with respect to the symplectic bases $\Sigma_B$ and
$\Sigma_L$, respectively. For an element $z\in \mathbf  C^6$, we define 
$z^{\#} = z(C\tau +D)^{-1}$.

For a rational vector $n=(n',n'') \in \mathbf  Q^{12}$,
we consider the following functional equations for a function
$F$ of $z^{\#}\in \mathbf  C^6$:
\begin{align*}
\text{ ($Eq_n^\#$): }
& F(z^{\#}+ p\tau^{\#} +q) =
\mathbf  e (-\frac{1}{2}p\tau^{\# t} p -p\;^tz^\#)
\mathbf  e(n^{'t}q - n^{''t}p)
F(z^{\#}), \\
\text{ ($Eq_m$): }
& F(z+ r\tau +s) =
\mathbf  e (-\frac{1}{2}r\tau r -r\;^tz)
\mathbf  e(m^{'t}s - m^{''t}r)
F(z).
\end{align*}
Note that the theta function $\vartheta_n(\Sigma_L,z^\#)$ satisfies
the functional equations \text{ ($Eq_n^\#$) }
for $p, q \in \mathbf  Z^6$.

Before studying the relation between $\vartheta_n(\tau^\#,z^\#)$
and $\vartheta_m(\tau,z)$, we define a translation vector
$\delta$ relative to the matrix $\sigma$.
We define $\delta' \in \mathbf  Z^6$ and 
$\delta'' \in \frac{1}{2}\mathbf  Z^3\oplus \mathbf  Z^3$ as
$$
\delta'=(\;^tC A)_0,\quad \delta''=(e\;^tD B e)_0 e^{-1}.
$$
The vector $(\delta',\delta'')$ is denoted by $\delta_{\Sigma}$.
To describe properties of $\delta_{\Sigma}$, it 
is convenient to consider a 
quadratic form $q$ on 
$\frac{1}{1-\rho}H_1(C, \mathbf  Z)^{-}/H_1(C, \mathbf  Z)^{-}$ induced
by the quadratic form $q$ on
$H_1(C, \mathbf  Z)^{-}/(1-\rho)H_1(C, \mathbf  Z)^{-}$  defined in 
\S \ref{total singular}
via the isomorphism
\begin{align*}
\frac{1}{1-\rho}H_1(C, \mathbf  Z)^{-}/H_1(C, \mathbf  Z)^{-}
& \to H_1(C, \mathbf  Z)^{-}/(1-\rho)H_1(C, \mathbf  Z)^{-} \\
x &  \mapsto (1-\rho)x.
\end{align*}

Let $(L_1,\Sigma_1)$ be the good basis defined in 
\S \ref{polarization of prym var} and
$g \in U(H_1(C, \mathbf  Z)^-)$. We put $L_g = g(L_1)$ and 
$\Sigma_g=g(\Sigma_1)$. Then it easy to see that $(L_g,\Sigma_g)$
is a good symplectic basis of $L_g$. The translation
vector $\delta_{\Sigma_g}$ is denoted by $\delta_g$.

\begin{proposition}
\begin{enumerate}
\item
The vector $\delta_g=(\delta', \delta'')$ is contained in
$2\mathbf  Z^3 \oplus \mathbf  Z^3 \oplus \mathbf  Z^6$.
\item
Let 
$\overline{\Delta}$ and 
$\overline{\frac{1}{2}\delta_g}$
be $e_1+e_2+e_5+e_6$ 
under the mapping (\ref{marking mod 2})
and the class
of  $\frac{1}{2}\delta_g\;^t(\alpha, \beta)$ 
in
$\frac{1}{1-\rho}H_1(C, \mathbf  Z)^{-}/H_1(C, \mathbf  Z)^{-}$.
Then
$$
c_0=
\overline{\frac{1}{2}\delta_g}-
\overline{\Delta} g
$$
is independent of $g$.
Moreover we have $q(\overline{\Delta})=0$ and $\overline{\Delta} \neq 0$. 
\end{enumerate}
\end{proposition}
\begin{proof}
For any $(L_g,\Sigma_g)$, we compute vectors $\delta_g$ by definition.  
As a consequence, we obtain this proposition.  
\end{proof}

For  $m=(m', m'') \in \mathbf  Q^{12}$, we set 
$\tilde m= m+\frac{1}{2}\delta_{\Sigma}, n=\tilde m\sigma^{-1}$. 
i.e., $n\;^t(\mathbf  a,\mathbf  b)=\tilde m\;^t(\alpha, \beta)$.
The next proposition is fundamental for comparing theta functions
for different lattices. We define $\theta_m(z^\#)$ by
$$
\theta_m(z^\#)=
\mathbf  e(\frac{1}{2}z(C\tau+D)^{-1}C\;^tz)\vartheta_m(\Sigma_B, z).
$$
In this definition, $z$ denotes the function of $z^\#$ 
by the relation $z^\# =z(C \tau +D)^{-1}$.

\begin{proposition}
\label{corresp. L and B}
\begin{enumerate}
\item
The function $\theta_m(z^{\#})$ satisfies the functional equation
($Eq_n^\#$) for $(p,q)\;^t(\mathbf  a, \mathbf  b) \in 
(1-\rho)H_1(C, \mathbf  Z)^-$. 
\item
The function $\theta_{m'}(z^\#)$ satisfies the functional equation
($Eq_{n}^\#$) for $m' \in \mathbf  Q^{12}$, 
$(p,q)\;^t(\mathbf  a, \mathbf  b) \in (1-\rho)H_1(C, \mathbf  Z)^-$
if and only if 
$m - m' \in \mathbf  Z^6 \oplus \frac{1}{2}\mathbf  Z^3
\oplus \mathbf  Z^3$, i.e.
$(m - m')\;^t(\alpha,\beta) \in H_1(C, \mathbf  Z)^-$. 
\end{enumerate}
\end{proposition}

Let $\Theta(\Sigma_B,m)$ and $\Theta(\Sigma_L,n)$ be the spaces of functions
of $\tau$ and $\tau^\#$ satisfying the
functional equations $(Eq_m)$ and $(Eq_n^\#)$
for $(r,s)\;^t(\alpha, \beta) \in (1-\rho)H_1(C, \mathbf  Z)^-$ and 
$(p,q)\;^t(\mathbf  a, \mathbf  b) \in (1-\rho)H_1(C, \mathbf  Z)^-$,
respectively. Since the space
$\Theta(\Sigma_B,m)$ and $\Theta(\Sigma_L,n)$ are $8$ dimensional
by the Riemann Roch theorem, 
Proposition \ref{corresp. L and B} implies 
the following proposition.
\begin{proposition}
\label{cor theta diff L}
Let $n=(m+\frac{1}{2}\delta_{\Sigma})\sigma^{-1}$.
By mapping a function $f(z)$ of $z$ to 
a function $f^\#(z^\#)=\mathbf  e(\frac{1}{2}z(C\tau+D)^{-1}C\;^tz)f( z)$,
of $z^\# = z(C \tau +D)^{-1}$, we have an isomorphism
$$
\Theta(\Sigma_B,m) \to \Theta(\Sigma_L,n).
$$
\end{proposition}

\subsection{$\Sigma$-trace}

In order to express
$\vartheta_0(\Sigma_{L},z^\#)$ as a linear
combination of translations of $\theta_m(z^\#)$ 
in the last section with
simple exponential coefficients,
we introduce $\Sigma$-trace.

\begin{definition}[$\Sigma$-trace]
Let $m$ be an element such that
$m+\frac{1}{2}\delta_{\Sigma} 
\in \mathbf  Z^6 \oplus \frac{1}{2}\mathbf  Z^3
\oplus \mathbf  Z^3$.
Let $S$ be a representative of $(p,q)$ for
$$
\mathbf  Z^{12}/\{(c,d) \mid (c,d)\;^t(\mathbf  a,\mathbf  b)\in 
(1-\rho)H_1(C, \mathbf  Z)^-\}.
$$
The $\Sigma$-trace $tr_\Sigma(f)(z^\#)$
is defined by
$$
tr_\Sigma(f)(z^\#)=\sum_{(p,q) \in S}
f(z^\# + p\tau^\# +q )
\mathbf  e (\frac{1}{2}p\tau^{\# t} p +p\;^tz^\#).
$$
\end{definition}

\begin{proposition}
The $\Sigma$-trace $tr_{\Sigma}(\theta_m)$
is independent of the choice of the representative $S$,
and it is a constant multiple of 
$\vartheta_0(\Sigma_L,z^\#)$.
Moreover there exists
$m \in \mathbf  Z^6 \oplus \frac{1}{2}\mathbf  Z^3
\oplus \mathbf  Z^3 -\frac{1}{2}\delta_{\Sigma}$ such that 
$tr_\Sigma(\theta_m)$ is non-zero.
\end{proposition}

\begin{proof}
Note that
the $\Sigma$-trace $tr_{\Sigma_L}(\theta_m)$ satisfies $(Eq_n^\#)$
for $(p,q)$ in a sufficiently small lattice in $L$.
By Proposition 
\ref{corresp. L and B}
and the characterization of the space generated by 
theta functions for principally polarized
abelian varieties, we have this proposition.
\end{proof}

\begin{definition}[$\Phi_g$,$\Phi_{g,n}$]
\label{def of phi}
\begin{enumerate}
\item
For each $g\in U(H_1(C, \mathbf  Z)^-)$, 
we choose $m_g
\in \mathbf  Z^6 \oplus \frac{1}{2}\mathbf  Z^3
\oplus \mathbf  Z^3-\frac{1}{2}\delta_g$
such that $tr_{\Sigma_g}(\theta_{m_g})$ is non-zero.
We define $\Phi_g$ by $tr_{\Sigma_g}(\theta_{m_g})$ and
$c_g =\displaystyle\frac{\vartheta_0(\Sigma_g,z^\#)}{\Phi_g(z^\#)}$.
\item
For $n= (p_0,q_0) \in \mathbf  Q^{12}$, 
we define 
$$
\Phi_{g,n}(z^\#)=
\mathbf  e(\frac{1}{2}p_0\tau \;^tp_0+p_0\;^t(z^\#+q_0))
\Phi_g(z^\#+p_0\tau^\# + q_0)
$$
and
$\Phi_{g,n}=\Phi_{g,n}(0)$.
\end{enumerate}
\end{definition}

Recall that $\Phi_{g, n}(z^\#)$ is a linear  
combination of translations of $\vartheta (\Sigma_B,z^\#)$.
In the rest of this subsection, we compute its coefficients.

\begin{proposition}
For $n_0=(m_0+\frac{1}{2}\delta_g)\sigma^{-1}$, 
we choose representatives $S_g(n_0)$ and $S_B(m_0)$ of $n_0+\{\tilde n \mid
\tilde n^{t}(\mathbf  a, \mathbf  b) 
\in H_1(C, \mathbf  Z)^- \}/\mathbf  Z^{12}$
and $m_0+\mathbf  Z^6 \oplus \frac{1}{2}\mathbf  Z^3 \oplus \mathbf  Z^3/\mathbf  Z^{12}$,
respectively. Then
$\{\Phi_{n}(z^\#)\}_{n \in S_g(n_0)}$
and $\{\vartheta_{m}(\Sigma_B,z)\}_{m \in S_B(m_0)}$
are bases of the $8$-dimensional vector spaces $\Theta (\Sigma_g,n_0)$
and $\Theta (\Sigma_B,m)$, respectively.
\end{proposition}

Let $(p_0,q_0) \in \mathbf  Q^{12}$. By the definition of 
$\Phi_{g,(p_0,q_0)}(z^\#)$, we have the following proposition 
by simple calculation.
\begin{proposition}
\label{L trace and theta const}
For $(p_0, q_0) \in \mathbf  Q^{12}$, we have
\begin{align*}
& \Phi_{g,(p_0,q_0)}(z^\#)  \\
= 
\sum_{(p,q)\in S} &
c_{(p_0,q_0),(p,q)}^{(g)}
\vartheta_{m_g+(p+p_0, q+q_0)\sigma_g}(\Sigma_B,z)
\cdot
\mathbf  e(\frac{1}{2}z(C\tau+D)^{-1}C\;^tz).
\end{align*}
and
\begin{align*}
\Phi_{g,(p_0,q_0)} =
\sum_{(p,q)\in S} 
c_{(p_0,q_0),(p,q)}^{(g)}
&\vartheta_{m_g+(p+p_0, q+q_0)\sigma_g}(\Sigma_B),
\end{align*}
where
$$
c_{(p_0,q_0),(p,q)}^{(g)} =
\mathbf  e(-\frac{1}{2}q\;^tp -\frac{1}{2}s\;^t r)
\mathbf  e(-r\;^tm_i'') 
\mathbf  e(\frac{1}{2}p_0\;^tq_0 -\frac{1}{2}r_0\;^ts_0 
-r_0\;^tm_i''-r_0\;^ts). 
$$
Here, $m_g$ and $S$ are in the definition of $\Sigma_g$-trace
and $\Phi_g(z^\#)$,
$(r_0,s_0)=(p_0,q_0)\sigma_g$ and
$(r,s)=(p,q)\sigma_g$.
\end{proposition}

Let $n=(p_0, q_0)$ and $(p,q)$ be elements of $S_g(n_0)$ and $\mathbf  Z^{12}$,
respectively. For
$(r_0,s_0)=(p_0, q_0) \sigma_g$ and $(r,s)=(p, q) \sigma_g$, 
we have
\begin{align*}
m_g+ (p+p_0, q+q_0)\sigma_g
\in & m_g+ n_0 \sigma_g +\mathbf  Z^6\oplus 
\frac{1}{2}\mathbf  Z^3 \oplus \mathbf  Z^3 \\
= & m_g+ m_0 + \frac{1}{2}\delta_g +\mathbf  Z^6\oplus 
\frac{1}{2}\mathbf  Z^3 \oplus \mathbf  Z^3 \\
= &  m_0 +\mathbf  Z^6\oplus 
\frac{1}{2}\mathbf  Z^3 \oplus \mathbf  Z^3. 
\end{align*}
For an element $m \in S_B(m_0)$, we put
$$ 
I(m)=\{(p,q)\in S\mid m_g+(p+p_0,q+q_0)\sigma - m \in \mathbf  Z^{12} \},
$$
$$
d_{n,m}^{(g,B)}=\sum_{(p,q) \in I(m)}c_{(p_0,q_0),(p,q)}^{(g)} 
\frac{\vartheta_{m_g+(p+q_0, q+q_0)\sigma}(\Sigma_B,z)}
{\vartheta_{m}(\Sigma_B,z)}.
$$
Then
by Proposition \ref{L trace and theta const}, for $n \in S_g(n_0)$,
$\Phi_{g,n}(z^\#)$ can be written as
\begin{equation}
\label{changing summation}
\Phi_{g,n}(z^\#) =
\sum_{m \in S_B(m_0)} 
d_{n,m}^{(g,B)}
\vartheta_{m}(\Sigma_B,z)
\cdot
\mathbf  e(\frac{1}{2}z(C\tau+D)^{-1}C\;^tz).
\end{equation}

We put 
$$
D^{(g,B)} = (d_{n,m}^{(g,B)})_{n \in S_g(n_0), m \in S_B(m_0)}.
$$
This is the base change matrix for
$\{\vartheta_m(\Sigma_B, z)\}_{m \in S_B(m_0)}$
and
$\{\Phi_{g,n}(z^\#)\}_{n \in S_g(n_0)}$
up to a constant exponential multiple.

\section{Main Theorem}
\subsection{Action of the stabilizer of length $0$ element
on theta functions}
\label{action of theta stabl.}

By choosing a principal lattice $L$, we get an injective homomorphism 
from $U(H_{std})$ to
$Sp(6,\mathbf  R) =Aut(\frak H_6)$ and
an inclusion $\jmath : B(H_{std}) \to \frak H_6$.
By this inclusion, we identify $U(H_{std})$ as a subgroup of
$Sp(6,\mathbf  R)$.
We consider the function $\det (C\tau +D)^{1/2}$ on 
$Sp(6, \mathbf  R)\times \frak H_6$ defined in 
\S \ref{application of quadr. rel of theta}.
Let $\sigma_g =\left(\begin{matrix} A_g & B_g \\ C_g & D_g \end{matrix}\right)$
be a matrix such that $\Sigma_{L_g}=\sigma_g(\Sigma_B)$.
We put $n_0=(m_0+\frac{1}{2}\delta_g)\sigma_g^{-1}$ for $m_0 \in \mathbf  Q^{12}$.
We choose $S_g(n_0)$ and $S_B(m_0)$ as in the last section.
For any $n \in S_L(n_0)$, there exist complex numbers
$u_{n,m}^{(g,B)}$ ($m \in S_B(m_0)$) independent of $\tau$ and $z$
such that
\begin{align}
\label{Igusa gransf. non-principal}
\vartheta_{n}(\Sigma_g,z^\#_g) = & \det(C_g \tau +D_g)^{-1/2}
\mathbf  e(z(C_g\tau +D_g)C_g^{-1 t}z) \\
& \sum_{m \in S_B(m_0)} u_{n, m}^{(g,B)}\vartheta_{m}(\Sigma_B,z),
\nonumber
\end{align}
where
$\tau^\#_g=(A_g\tau +B_g)(C_g\tau+D_g)^{-1}$ and 
$z^\#_g=z\cdot (C_g\tau +D_g)^{-1}$
by the transformation formula in p.84, \cite{I}.
Moreover this expression is unique.
By comparing the right hand sides of (\ref{changing summation}) and
(\ref{Igusa gransf. non-principal}), we have the
following proposition.
\begin{proposition}
\label{homothety}
The matrix $U^{(g,B)}=(u_{n,m}^{(g,B)})_{n \in S_L(n_0),m \in S_B(m_0)}$
is a non-zero constant multiple $c_g$ of 
$D^{(g,B)}$.
Moreover the non-zero constant $c_g$ does not depend on $\tau$.
Especially for $n \in S_g(0)$, we have
\begin{align}
\label{theta function transform dif lat}
\vartheta_{n}(\Sigma_g,z^\#_g) & = \sum_{m \in S_1(0)}
u^{(g)}_{n,m}\vartheta_{m}(\Sigma_1,z^\#_1) \\
& = c_g\sum_{m \in S_1(0)}
d^{(g)}_{n,m}\vartheta_{m}(\Sigma_1,z^\#_1). \nonumber
\end{align}

\end{proposition}
We define
$$
U(H_{std})_{\overline{\Delta}}=\{ g \in U(H_{std}) \mid
\overline{\Delta}g =\overline{\Delta}\}.
$$
We fix a representative $S_1(0)$.
For example we choose 
$S_1(0)=\{\frac{1}{2}(\mu_j, \mu_jU)\}$ with
{\small
\begin{align*}
&
\mu_1=(0,0,0,0,0,0),
\mu_2=(0,0,1,1,1,1),
\mu_3=(1,1,0,0,1,1),
\mu_4=(1,1,1,1,0,0), \\
&
\mu_5=(1,1,1,1,1,1),
\mu_6=(1,1,0,0,0,0),
\mu_7=(0,0,1,1,0,0), 
\mu_8=(0,0,0,0,1,1).
\end{align*}
}
If $n\;^t(\mathbf  a, \mathbf  b) \in H_1(C, \mathbf  Z)^-$, 
then $n\;^t(g(\mathbf  a), g(\mathbf  b))=
g(n\;^t(\mathbf  a, \mathbf  b)) \in g(H_1(C, \mathbf  Z)^-)$ for 
$g \in U(H_{std})_{\overline{\Delta}}$.
Therefore we can take a representative $S_g(0)$ as $S_1(0)$.
Then the vector spaces generated by 
$\Phi_{n}(z^\#_g)$ ($n \in S_g(0)$) and 
$\Phi_{n}(z^\#)$ ($n \in S_1(0)$)
are isomorphic via the map 
defined in Proposition \ref{cor theta diff L}:
$$
\Theta(\Sigma_g,0)\overset{\simeq}{\longleftarrow}
\Theta(\Sigma_B,-\overline{\Delta}) \overset{\simeq}{\longrightarrow}
\Theta(\Sigma_1,0), 
$$
where $z^\#_g$ and $z^\#$ are related by
$z^\# = z(C\tau +D)^{-1}$ and 
$z^\#_g = z(C_g\tau +D_g)^{-1}$.
We put $D^{(g)}=D^{(g,B)}(D^{(id,B)})^{-1}$ and
$U^{(g)}=U^{(g,B)}(U^{(id,B)})^{-1}$.
By the definition of $U^{(g,B)}$, the map
$g \in U(H_{std})_{\overline{\Delta}} \mapsto U^{(g)}$
defines a projective representation of $U(H_{std})_{\overline{\Delta}}$,
which is denoted by $\chi$.

Since $U^{(g)}=c_g c_1^{-1}D^{(g)}$,
we have the following corollary of Proposition \ref{homothety}.

\begin{corollary}
The map 
$$
U(H_{std})_{\overline{\Delta}} \ni g \mapsto
D^{(g,B)}(D^{(1,B)})^{-1} \in Aut(\bold C^{S_1(0)})
$$ 
becomes a projective representation of
$U(H_{std})_{\overline{\Delta}}$, which is isomorphic to $\chi$.
\end{corollary}

Let $S_1(0)_{ev}$ be the subset of $S_1(0)$ consisting of
$v\;^tUv\in 4\mathbf  Z$. 
In the example given as above, we have  
$S_1(0)_{ev}=\{\frac{1}{2}(\mu_j, \mu_jU)\}$ with $j=1,4,6,7$.
By evaluating 
\ref{theta function transform dif lat} at $z^\#_1=z^\#_g=0$,
for $n_2 \in S_1(0)_{ev}$, we have
\begin{align}
\label{theta const transform dif lat}
\vartheta_{n}(\Sigma_g) & = \sum_{m \in S_1(0)_{ev}}
u^{(g)}_{n,m}\vartheta_{m}(\Sigma_1) \\
& = c_gc_1^{-1}\sum_{m \in S_1(0)_{ev}}
d^{(g)}_{n,m}\vartheta_{m}(\Sigma_1). \nonumber
\end{align}
We put $D^{(g)}_{ev}=(d^{(g,B)}_{n, m})_{n,m \in S_1(0)_{ev}}$.
We define a projective representation $\chi_{const}$ as 
$$
\chi_{const}(g)=D^{(g)}_{ev}
 \in PGL(4, \mathbf  C)
$$
on the space of theta constants.
Note that an element $g$ in $U(H_{std})$ is in 
$U(H_{std})_{\overline{\Delta}}$ if and only if 
$\pi (g) \in (\frak S_4(1,2,5,6) \times \frak S_4(3,4,7,8))\rtimes \frak S_2$
under the homomorphism
$\pi :U(H_{std}) \to \frak S_8$
defined just after Lemma \ref{O6+(2) to S8}.
Here $\frak S_4(1,2,5,6)$ is 
the symmetric group of permutations of index $\{1,2,5,6\}$.
Let $M_{2,5}$ be the (complex) reflection corresponding to
the transposition of the points $p_2$ and $p_5$.
Then we have $M_{2,5} \in U(H_{std})_{\overline{\Delta}}$,
$$
D^{(M_{2,5})}_{ev}=
\left(\begin{matrix}
\frac{1}{2}-\frac{1}{2}i & - \frac{1}{2}-\frac{1}{2}i & 0 & 0 \\
-\frac{1}{2}-\frac{1}{2}i &  \frac{1}{2}-\frac{1}{2}i & 0 & 0 \\
0 & 0 & \frac{1}{2}-\frac{1}{2}i & - \frac{1}{2}-\frac{1}{2}i  \\
0 & 0 & -\frac{1}{2}-\frac{1}{2}i &  \frac{1}{2}-\frac{1}{2}i  
\end{matrix}\right) 
$$
and
\begin{align*}
 \;^t(\vartheta_{m_j}(\tau^\#_g))_{j=1,4,3,2}  
= & 
c_{M_{2,5}}c_1^{-1}\cdot
\det (\gamma_g\tau^\# + \delta_g)^{1/2} \\
& \cdot D_{ev}^{(M_{2,5})t}
(\vartheta_{m_j}(\tau^\#_1))_{j=1,4,3,2},
\end{align*}
where 
$$
\sigma_g\cdot \sigma^{-1}_1 =
\left(\begin{matrix}
\alpha_g & \beta_g \\ \gamma_g & \delta_g 
\end{matrix}\right) \in Sp(6, \mathbf  Q).
$$

\subsection{Theta constants and cross-ratios of coordinates}

In this section, we combine the results in 
\S \ref{determination} and
\S \ref{action of theta stabl.}
to get an $\frak S_8$-equivariant presentation of a projective
map from the moduli space $M_{8pt}$ to $\mathbf  P^{104}$.

We define
a function $\Cal T_g(\tau)$ of $\tau \in M_{marked}$ by 
\begin{equation}
\label{def tg}
\Cal T_g(\tau)=
\det(\gamma_g \tau +\delta_g)^{-1}
 \vartheta_{m_1}(\tau^\#_g)
\vartheta_{m_3}(\tau^\#_g),
\end{equation}
where $\sigma_g=\left(\begin{matrix}
\alpha_g& \beta_g\\
\gamma_g& \delta_g\\
\end{matrix}\right)$ 
and $\tau^\#_g$ are defined as before for $g\in U(H_{std})$.
Then we have
\begin{equation}
\label{equivariance for tg}
\Cal T_g(h\cdot\tau)=\Cal T_{gh}(\tau).
\end{equation}
Since $\Cal T_g^2$ depends only on the image $\pi (g) \in \frak S_8$,
it is also denoted by $\Cal T_{\pi (g)}^2$.

By the result of the last section, if 
$g \in U(H_{std})_{\overline{\Delta}}$, then $\Cal T_g(\tau)$
is a homogeneous polynomial of $\vartheta_{m_j}(\tau)$
with constant coefficients.
For example if $g= M_{25}$, 
we have
\begin{align*}
\Cal T_g =  \frac{-i}{2}
\cdot (c_{M_{2,5}}c_1^{-1})^2
& (\vartheta_{m_1}(\tau^\#_1)-
i\vartheta_{m_4}(\tau^\#_1)) \\
& \cdot(\vartheta_{m_3}(\tau^\#_1)
-i\vartheta_{m_2}(\tau^\#_1)).
\end{align*}
Therefore we have
\begin{align}
\label{equality for quot. square of T}
 \frac{\Cal T_g^2}{\Cal T_1^2} 
= &  \frac{-1}{4}
\cdot (c_{M_{2,5}}c_1^{-1})^4
\frac{
(\vartheta_{m_1}(\tau^\#_1)-
i\vartheta_{m_4}(\tau^\#_1))^2}
{\vartheta_{m_1}(\tau^\#_1)^2} 
 \\
&\cdot\frac
{(\vartheta_{m_3}(\tau^\#_1)
-i\vartheta_{m_2}(\tau^\#_1))^2
}
{\vartheta_{m_3}(\tau^\#_1)^2} 
\nonumber \\
= & c\frac{(x_1-x_5)(x_2-x_6)}{(x_1-x_2)(x_5-x_6)},
\nonumber
\end{align}
where $c=(c_{M_{2,5}}c_1^{-1})^4$.

Let $R$ be a set of representatives 
of the composite surjection
$$
U(H_{std}) \overset{\pi}\longrightarrow \frak S_8 \to
Stab\{\{1,2\},\{5,6\},\{3,4\},\{7,8\}\} \backslash \frak S_8,
$$
where $\pi$
is the natural surjection and
$Stab\{\{1,2\},\{5,6\},\{3,4\},\{7,8\}\}$ is 
the stabilizer of $\{\{1,2\},\{5,6\},\{3,4\},\{7,8\}\}$.  
We fix this set $R$ once and for all.
\begin{definition}[Polynomial map]
\label{def polynomial map}
Set $P_1=(x_1-x_2)(x_3-x_4)(x_5-x_6)(x_7-x_8)$ and
$P_r=\pi(r)^*(P_1)$ for $r \in R$.
Since each $P_r$ is relative invariant under the action of $PGL(2,\mathbf C)$,
the map $P:(\mathbf  P^1)^8-\Diag \to \mathbf  P^{104}$ defined 
by the ratio of $(P_r)_{r \in R}$ descends to 
a morphism $M_{8pts} \to \mathbf  P^{104},$
which is also denoted by $P$. The composite $M_{marked} \to M_{8pts}
\to \mathbf  P^{104}$ is also denoted by $P$.
\end{definition}

By the Definition \ref{def polynomial map}, the last term of
(\ref{equality for quot. square of T}) is equal to
$c\cdot\frac{P_g}{P_1}$.

We have the following theorem. 
\begin{theorem}
\label{main I}
Let $\Cal T^{(2)}$ be the map from $M_{marked}$ to 
$\mathbf  P^{104}$ defined by $(\Cal T_r^2)_{r \in R}$.
Then the following diagram is commutative:

\setlength{\unitlength}{0.75mm}
\begin{picture}(200,45)(-40,-5)
\put(0,0){$B(H_{std})$}
\put(38,0){$B(H_{std})/\Gamma (1+i)$}
\put(0,30){$M_{marked}$}
\put(38,30){$M_{8pts}$}
\put(71,14.5){$\mathbf  P^{104}$.}
\put(6,15){$p$}
\put(55,27){$P$}
\put(30,11){$\Cal T^{(2)}$}
\put(3,25){\vector(0,-1){20}}
\put(42,25){\vector(0,-1){20}}
\put(18,32){\vector(1,0){18}}
\put(18,2){\vector(1,0){18}}
\put(48,4){\vector(2,1){20}}
\put(48,28){\vector(2,-1){20}}
\put(17,4){\vector(4,1){50}}
\end{picture}
\end{theorem}

In order to prove Theorem \ref{main I}, we give some lemmas. 
Let $[r]$ denotes the class of $r \in \frak S_8$ in
$Stab\{\{1,2\},\{5,6\},\{3,4\},\{7,8\}\} \backslash \frak S_8$.

\begin{lemma}
\begin{enumerate}
\item
If $[r]=[r']$, then 
$\Cal T_{r}^2$ is a constant multiple of $\Cal T_{r'}^2$.
\item
The map from $Stab\{\{1,2\},\{5,6\},\{3,4\},\{7,8\}\}$
to $\{ \pm 1\}$ defined by 
$g \mapsto \frac{\Cal T_g^2}{\Cal T_1^2}$
is a character  and
this coincides with the restriction of the signature on
$\frak S_8$.
\end{enumerate}
\end{lemma}
\begin{proof}
1. By the equality \ref{equivariance for tg},
We have the following equation of rational functions of
$M_{8pts}$,  
\begin{equation}
\label{quot of theta is a rational function}
\frac{\Cal T_g^2}{\Cal T_1^2}(h(\tau))=
\frac{\Cal T_{gh}^2}{\Cal T_h^2}(\tau)
\end{equation}
for $h \in \frak S_8$. 
Thus we have only to prove the lemma for the case $r'=1$.
If $[r]=[1]$, $\Cal T_{r}^2$ is a constant multiple of $\Cal T_1^2$
by Corollary \ref{cor to quadratic rel final}
and the expression of the projective
representation $\chi_{const}$.

2.  
The first statement is a consequence of 1.
Using the transformation formula of p.85 in \cite{I},
we have
\begin{equation}
\label{first character property}
\Cal T_{M_{12}}=\Cal T_{M_{56}}=i\Cal T_1.
\end{equation}
By applying $M_{12}M_{56}$ to the equality
(\ref{equality for quot. square of T}),
we have
\begin{equation}
\label{second character property}
\Cal T_{gM_{12}M_{56}}^2=\Cal T_{g}^2.
\end{equation}
Equalities (\ref{first character property}) and
(\ref{second character property}) characterize the character of
$Stab\{\{1,2\},\{5,6\},\{3,4\},\{7,8\}\}$ and the
character $\displaystyle \frac{\Cal T_{g}^2}{\Cal T_1^2}$
coincides with the restriction of the signature.
\end{proof}

\begin{lemma}
\label{chain}
Let $[r]$ be an element of
$Stab\{\{1,2\},\{5,6\},\{3,4\},\{7,8\}\}\backslash \frak S_8$ 
and $(2,6) \in \frak S_8$ be the transposition of $2$ and $6$.
Then there exist sequences of
$g_1, \dots ,g_{k+1}$ and $h_1, \dots ,h_k$ of $\frak S_8$ such that
\begin{enumerate}
\item
$[r]=[g_1h_1\cdots g_kh_kg_{k+1}]$,
\item
$[g_1h_1\cdots g_l]=[g_1h_1\cdots g_lh_l]$
for $l=1, \dots, k$,
\item
$[(2,6)]=[(2,6)g_1]$,
$[(2,6)g_1h_1\cdots g_lh_l]
=[(2,6)g_1h_1\cdots g_lh_lg_{l+1}]$ for $l=1, \dots, k$.
\end{enumerate}
\end{lemma}
\begin{proof}[Proof of Theorem \ref{main I}]
Using (\ref{quot of theta is a rational function}), we have
\begin{equation}
\label{comp 2,5 and 2,6 I}
\frac{\Cal T_{M_{2,5}}^2}{\Cal T_1^2}(M_{2,6}(\tau))=
\frac{\Cal T_{M_{2,5}M_{2,6}}^2}{\Cal T_{M_{2,6}}^2}(\tau)=
-\frac{\Cal T_{M_{2,5}}^2}{\Cal T_{M_{2,6}}^2}(\tau). 
\end{equation}
We put $c=(c_{M_{2,5}}c_1^{-1})^4$.
Since the map $p$ is equivariant under the action of $\frak S_8$, 
we have
\begin{equation}
\label{comp 2,5 and 2,6 II}
\frac{\Cal T_{M_{2,5}}^2}{\Cal T_1^2}(h(\tau))=
c\cdot\frac{(x_{1h}-x_{5h})(x_{2h}-x_{6h})}
{(x_{1h}-x_{2h})(x_{5h}-x_{6h})}
\end{equation}
for $h \in \frak S_8$. 
The equations (\ref{comp 2,5 and 2,6 I}) and
(\ref{comp 2,5 and 2,6 II}) yield 
\begin{equation}
\label{seeds of equation}
\displaystyle
\frac{\Cal T_{(2,6)}^2}{\Cal T_1^2}(\tau) =
\frac{\Cal T_{M_{2,6}}^2}{\Cal T_1^2}(\tau) =
\frac{P_{(2,6)}}{P_{1}}.
\end{equation}
For any $r\in frak S_8$, there exist sequences $g_1, \dots, g_{k+1},
h_1,\dots ,h_k$ such that
$$
\begin{matrix}
[(2,6)]&=&[(2,6)g_1],& & [(2,6)g_1h_1] &=& [(2,6)g_1h_1g_2], \\ 
& & [g_1] &=& [g_1h_1],& &[g_1h_1g_2] &=& [g_1h_1g_2h_2] 
\end{matrix}
$$
\vskip 0.1in
$$
\begin{matrix}
& & & \dots & [(2,6)g_1\cdots h_k] &=& [(2,6)g_1\cdots h_kg_{k+1}] \\ 
& &  & & & \dots &[g_1h_1 \cdots h_k g_{k+1}] &=& [r]. 
\end{matrix}
$$
by Lemma \ref{chain}. By applying $g_1$ and $g_1h_1$ to the equality
(\ref{seeds of equation}), we have
\begin{equation}
\label{seeds2 of equation}
\frac{\Cal T_{(2,6)}^2}{\Cal T_{g_1}^2}(\tau) =
\frac{\Cal T_{(2,6)g_1}^2}{\Cal T_{g_1}^2}(\tau) =
\frac{P_{(2,6)g_1}}{P_{g_1}}=
\frac{P_{(2,6)}}{P_{g_1}},
\end{equation}
\begin{equation}
\label{seeds3 of equation}
\frac{\Cal T_{(2,6)g_1h_1}^2}{\Cal T_{g_1}^2}(\tau) =
\frac{\Cal T_{(2,6)g_1h_1}^2}{\Cal T_{g_1h_1}^2}(\tau) =
\frac{P_{(2,6)g_1h_1}}{P_{g_1h_1}}=
\frac{P_{(2,6)g_1h_1}}{P_{g_1}}.
\end{equation}
From the equalities
(\ref{seeds of equation}) and
(\ref{seeds2 of equation}), we have
\begin{equation}
\label{seeds4 of equation}
\frac{\Cal T_{g_1}^2}{\Cal T_{1}^2}(\tau) =
\frac{P_{g_1}}{P_{1}}.
\end{equation}
and from the equalities
(\ref{seeds3 of equation}) and
(\ref{seeds4 of equation}), we have
$$
\frac{\Cal T_{(2,6)g_1h_1}^2}{\Cal T_{1}^2}(\tau) =
\frac{P_{(2,6)g_1h_1}}{P_{1}}.
$$
We continue this procedure, we get an identity
$
\displaystyle
\frac{\Cal T_{g}^2}{\Cal T_1^2}(\tau) =
\frac{P_{g}}{P_{1}}
$
for all $g \in \frak S_8$,
which completes the proof.
\end{proof}

\subsection{Branched covering of $M_{8pts}$ corresponding to
$\Gamma (2)$}

In this section, we study the map from $M_{marked}$ to
$\tilde M_{8pts}$ defined by the theta constants on $B(H_{std})$.
As in \S \ref{level2 str and 2 expo. cov},
we choose an initial point $X=(x_1, \dots ,x_8)$
and specify the branch of the function $\sqrt{x_k-x_j}$.

Let $r_1$ be $2^4$-partition $\{\{1,2,\}, \{3,4\},\{5,6\},\{7,8\}\}$
and $\pi : U(H_{std}) \to \frak S_8$ be the natural projection.
Using the argument $arg(g)$ defined in \S \ref{level2 str and 2 expo. cov},
we define  a multi-valued function $Q_r$ 
($r \in R$),
on $\mathbf  C^8-\Diag$ by
$$
Q_r=arg(g)\sqrt{(x_{j_2}-x_{j_1})(x_{j_4}-x_{j_3})
(x_{j_6}-x_{j_5})(x_{j_8}-x_{j_7})},
$$
where $r_1\pi(r)=
\{\{j_1,j_2\},\{j_3,j_4\},\{j_5,j_6\},\{j_7,j_8\}\}$
and $j_p < j_{p+1}$ for $j=1,3,5,7$.
Here we chose the branch of the square root as in 
\S \ref{level2 str and 2 expo. cov}. 
Let $\tilde N$ be the covering of $\mathbf  C^8-\Diag$
defined by $\sqrt{x_j-x_k}$ ($1 \leq k <j \leq 8$).
Then the functions $Q_r$ on $\tilde N$ define a morphism
$Q = (Q_r)_{r \in R}: \tilde N \to \mathbf  P^{104}$. 
Let $pr_1 :\mathbf  P^{104} \to \mathbf  P^{104}$ be the morphism
defined by $(y_r)_{r \in R} \mapsto (y^2_r)_{r \in R}$.
Since the coordinates of the inverse image of 
$P(\lambda_4, \dots, \lambda_8)$ under $pr_1$
can be expressed by polynomials of
$\sqrt{\lambda_j}, \sqrt{1-\lambda_j}, \sqrt{\lambda_j - \lambda_k}$,
the morphism $\tilde M_{8pts} \to M_{8pts} \to \bold P^{104}$
factors through $pr_1$ (see the following diagram).

\setlength{\unitlength}{0.75mm}
\begin{picture}(200,45)(-20,-5)
\put(7,30){$\tilde N$}
\put(15,32){\vector(1,0){21}}
\put(40,30){$\tilde M_{8pts}$} \put(58,32){\vector(1,0){18}}
  \put(80,30){$ \mathbf  P^{104}$}
\put(15,25){\vector(1,-1){20}} \put(64,26){$Q$}
\put(43,25){\vector(0,-1){20}}  \put(82,25){\vector(0,-1){20}}
\put(64,5){$P$}\put(87,15){$pr_1$}
\put(40,0){$M_{8pts}$}  \put(58,2){\vector(1,0){18}}
  \put(80,0){$\mathbf  P^{104}$.}
\end{picture}
We define the morphism $Q$ by the above diagram.

\begin{theorem}
Let $\Cal T=(\Cal T_{r})_{r \in R}$ be the morphism from
$B(H_{std})$ to $\mathbf  P^{104}$ defined by $\Cal T_r$ 
($r \in R$). 
Then the following diagram is commutative:

\setlength{\unitlength}{0.75mm}
\begin{picture}(200,45)(-40,-5)
\put(0,0){$B(H_{std})$}
\put(38,0){$B(H_{std})/\Gamma (2)$}
\put(0,30){$M_{marked}$}
\put(38,30){$\tilde M_{8pts}$}
\put(71,14.5){$\mathbf P^{104}$.}
\put(6,15){$p$}
\put(55,27){$Q$}
\put(30,11){$\Cal T$}
\put(3,25){\vector(0,-1){20}}
\put(42,25){\vector(0,-1){20}}
\put(18,32){\vector(1,0){18}}
\put(18,2){\vector(1,0){18}}
\put(48,4){\vector(2,1){20}}
\put(48,28){\vector(2,-1){20}}
\put(17,4){\vector(4,1){50}}
\end{picture}

\end{theorem}
\begin{proof} 
We determine the branch of the square root of the last term 
in (\ref{equality for quot. square of T}).
Let $h=M_{1,5},g=M_{2,5}$ be elements of $U(H_{std})$. 
Then we have $h(L_1)=gh(L_1)$.
By the transformation formula p.85 \cite{I}, we have 
$$
\Cal T_{gh} = i \Cal T_{h}.
$$ 
By Theorem \ref{main I}, we have
\begin{equation}
\label{first seeds}
\frac{\Cal T_{h}}{\Cal T_1}=i\cdot\epsilon\cdot
\sqrt{\frac{(x_5-x_2)(x_6-x_1)}{(x_2-x_1)(x_6-x_5)}}
\end{equation}
with $\epsilon =\pm 1$. By applying $g^*$ to the both sides of 
(\ref{first seeds}), and we get
\begin{equation}
\label{second seeds}
i\cdot \frac{\Cal T_{h}}{\Cal T_g}=
\frac{\Cal T_{hg}}{\Cal T_g}=
\frac{-\epsilon}{arg(g)}\cdot
\sqrt{\frac{(x_5-x_2)(x_6-x_1)}{(x_5-x_1)(x_6-x_2)}}.
\end{equation}
By (\ref{first seeds}) and (\ref{second seeds}), we have
\begin{equation*}
\frac{\Cal T_{g}}{\Cal T_1}=
\frac{\Cal T_{g}}{\Cal T_{h}}\cdot
\frac{\Cal T_{h}}{\Cal T_1}=
arg(g)\cdot
\sqrt{\frac{(x_5-x_1)(x_6-x_2)}{(x_2-x_1)(x_6-x_5)}}.
\end{equation*}
Using Lemma \ref{chain} and the same argument in the proof of
Theorem \ref{main I}, we have the theorem.
\end{proof}

\end{document}